\documentclass[journal]{IEEEtran}

\usepackage{cite}
\usepackage{amsfonts}
\usepackage{amsmath}
\usepackage{graphicx}
\usepackage{booktabs}
\usepackage{mathtools}
\usepackage{enumerate} 
\usepackage{multirow} 
\usepackage{mathrsfs}
\usepackage{ulem}
\usepackage{setspace}
\usepackage[ruled, linesnumbered]{algorithm2e}
\usepackage{enumitem}
\usepackage{tabulary}
\usepackage{tabularx}
\usepackage[dvipsnames]{xcolor}
\allowdisplaybreaks[4]
\usepackage{bm}
\usepackage{mathtools}
\usepackage{romannum}
\usepackage{url}
\usepackage{array}
\newcolumntype{L}{>{\centering\arraybackslash}m{2cm}}

\let\oldnl\nl
\newcommand{\nonl}{\renewcommand{\nl}{\let\nl\oldnl}}

\makeatletter
\newcommand{\removelatexerror}{\let\@latex@error\@gobble}
\makeatother

\begin{document}
\pagenumbering{arabic}	
\setlength{\textfloatsep}{3pt} 
\setlength{\abovecaptionskip}{3pt}
\setlength{\belowcaptionskip}{10pt}

\setlength{\abovedisplayskip}{-8.5pt}
\setlength{\belowdisplayskip}{4.5pt}

\title{Rolling Optimization of Mobile Energy Storage Fleets for Resilient Service Restoration}
\author{
	Shuhan~Yao,~\IEEEmembership{Student Member,~IEEE,}
	Peng~Wang,~\IEEEmembership{Fellow,~IEEE,}
	Xiaochuan~Liu,~\IEEEmembership{Student Member,~IEEE,}
	Huajun~Zhang,~\IEEEmembership{Student Member,~IEEE,}
	and~Tianyang~Zhao,~\IEEEmembership{Member,~IEEE}
	
\thanks{This work was supported by the Future Resilient Systems (FRS) at the Singapore-ETH Centre (SEC), which is funded by the National Research Foundation of Singapore (NRF) under its Campus for Research Excellence and Technological Enterprise (CREATE) program. \textit{(Corresponding author: Tianyang Zhao)} }
\thanks{S. Yao and H. Zhang are with the Interdisciplinary Graduate School, Nanyang Technological University, Singapore 639798 (email: syao002@e.ntu.edu.sg; hzhang031@e.ntu.edu.sg)}
\thanks{P. Wang and X. Liu are with the School of Electrical and Electronic Engineering, Nanyang Technological University, Singapore 639798 (e-mail: epwang@ntu.edu.sg; e160064@e.ntu.edu.sg).}
\thanks{T. Zhao is with the Energy Research Institute@NTU, Nanyang Technological University, Singapore 639798 (e-mail: zhaoty@ntu.edu.sg).} 
}

%
%
\markboth{IEEE TRANSACTIONS ON SMART GRID}
{Shell \MakeLowercase{\textit{et al.}}: Bare Demo of IEEEtran.cls for IEEE Journals}

\maketitle

\begin{abstract}
Mobile energy storage systems (MESSs) provide promising solutions to enhance distribution system resilience in terms of mobility and flexibility. 
This paper proposes a rolling integrated service restoration strategy to minimize the total system cost by coordinating the scheduling of MESS fleets, resource dispatching of microgrids and network reconfiguration of distribution systems. 
The integrated strategy takes into account damage and repair to both the roads in transportation networks and the branches in distribution systems. 
The uncertainties in load consumption and the status of roads and branches are modeled as scenario trees using Monte Carlo simulation method. 
The operation strategy of MESSs is modeled by a stochastic multi-layer time-space network technique. A rolling optimization framework is adopted to dynamically update system damage, and the coordinated scheduling at each time interval over the prediction horizon is formulated as a two-stage stochastic mixed-integer linear program with temporal-spatial and operation constraints.
The proposed model is verified on two integrated test systems, one is with Sioux Falls transportation network and four 33-bus distribution systems, and the other is the Singapore transportation network-based test system connecting six 33-bus distribution systems. 
The results demonstrate the effectiveness of MESS mobility to enhance distribution system resilience due to the coordination of mobile and stationary resources.


\end{abstract}


\IEEEpeerreviewmaketitle
\section*{Nomenclature}
\addcontentsline{toc}{section}{Nomenclature}
\subsection*{Sets}
\begin{IEEEdescription}[\IEEEusemathlabelsep\IEEEsetlabelwidth{$~~~~~~~~~~~~~~~~~$}]
	\item[$\mathcal{T}$] Set of time intervals. 
	\item[$\mathcal{S}$] Set of scenarios.
	\item[$\mathcal{N}_\text{D}$] Set of buses in distribution systems.
	\item[$\mathcal{E}_\text{D}$] Set of branches in distribution systems.
	\item[$\ddot{\mathcal{E}}_\text{D}^s$] Set of damaged branches in scenario $s$.
	\item[$\mathcal{N}_\text{T}$] Set of nodes in transportation networks.
	\item[$\mathcal{E}_\text{T}$] Set of edges in transportation networks. 
	\item[$\ddot{\mathcal{E}}_\text{T}^s$] Set of damaged edges in scenario $s$.
	\item[$\mathcal{M}$] Set of microgrids.
	\item[$\mathcal{D}$] Set of depots.
	\item[$\Omega$] Set of MESS fleet.
	\item[$\mathcal{N}_\text{S}^{\omega, s}$] Set of nodes in time-space networks for MESS $\omega$ in scenario $s$.
	\item[$\mathcal{E}_\text{S}^{\omega, s}$] Set of arcs in time-space networks for MESS $\omega$ in scenario $s$.
	
\end{IEEEdescription}

\subsection*{Parameters}	
\begin{IEEEdescription}[\IEEEusemathlabelsep\IEEEsetlabelwidth{$~~~~~~~~~~~~~~~~~$}]
	\item[$\Delta t$] Length of time intervals.
	\item[$r_{ij}$, $x_{ij}$] Resistance and reactance of line $(i,j)$.
	\item[$V_{\text{avg}}^\omega$] MESS $\omega$'s average speed.
	\item[$\varphi_i$] Power factor of load at bus $i$.
	\item[$\overline{P}_\text{ch}^\omega$ and $\overline{P}_\text{dch}^\omega$] Maximum charging/discharging power of MESSs.
	\item[ $\eta_\text{ch}^\omega$, $\eta_\text{dch}^\omega$] Charging/discharging efficiency.
	\item[ $E_\text{c}^\omega $]  Battery's capacity of MESSs.
	\item[$\overline{SOC}^\omega$, $\uline{SOC}^\omega$] Maximum and minimum level of SOC of the MESS.
	\item[$\overline{S}_{ij}$] Branch power capacity.
	\item[$\overline{v}_i$, $\uline{v}_i$] Maximum/minimum voltage magnitude.
	\item[$\varphi_i$] Power factor.
	\item[$\gamma_s$] Probability of scenario $s$.
	\item[$P_{\text{D}, i}^{t, s}$, $Q_{\text{D}, i}^{t, s}$] Predicted value of active/reactive load at bus $i$ in interval $t$ in scenario $s$.
	\item[$P_{\text{D}, i}^{t, \xi}$, $Q_{\text{D}, i}^{t, \xi}$] Realization of active/reactive load at bus $i$ in interval $t$.
	\item[$\overline{P}_{\text{DG},m}$, $\overline{Q}_{\text{DG},m}$] Maximum active/reactive power of equivalent dispatchable DG.
	\item[$\overline{E}_{\text{DG}, m}$, $\uline{E}_{\text{DG}, m}$] Energy capacity and minimum reserve in the microgrid. 
	\item[ $W_i$]  Unit interruption cost for load at bus $i$.
	\item[$C_{\text{gen}, m}$] Unit generation cost of microgrid $m$.
	\item[$C_{\text{bat}, \omega}$] Unit battery maintenance cost for MESS $\omega$.
	\item[$C_{\text{tran}, \omega}$] Unit transportation cost for MESS $\omega$.
\end{IEEEdescription}

\subsection*{Variables}
\begin{IEEEdescription}[\IEEEusemathlabelsep\IEEEsetlabelwidth{$~~~~~~~~~~~~~~~~~$}]
	\item[$P_{\text{ch},m}^{\omega, t, s}, P_{\text{dch},m}^{\omega, t, s}$] Charging/discharging power of MESS $\omega$ from/to microgrid $m$ in interval $t$ in scenario $s$.
	\item[$I_\text{ch}^{\omega, t, s}$, $I_\text{dch}^{\omega, t, s}$ ] Binary variables, battery status of MESS $\omega$ in interval $t$ in scenario $s$, 1 if the status is on, 0 otherwise.
	\item[$E_\omega^{t, s}$] Energy of the MESS $\omega$ by the end point of interval $t$ in scenario $s$.
	\item[$\alpha_{ij}^s$ ] Binary variable, the connection status of branch $(i, j)$ in scenario $s$, 1 is connected, 0 otherwise.
	\item [$\zeta _{\hat{n}\check{n}}^{\omega, s}$] Binary variables, 1 if the MESS $\omega$  is on an arc $(\hat{n}, \check{n}) \in \mathcal{E}_\text{S}^{\omega, s}$ in scenario $s$, otherwise set to 0. 
	\item[$P_{\text{G}, i}^{t, s}, Q_{\text{G}, i}^{t, s}$] Real/reactive power generation at bus $i$ in interval $t$ in scenario $s$.
	\item[$P_{\text{r}, i}^{t, s}$, $Q_{\text{r}, i}^{t, s}$] Load restoration at bus $i$ in interval $t$ in scenario $s$.
	\item[$P_{ij}^{t, s}$, $Q_{ij}^{t, s}$] Real/reactive power of branch $(i,j)$ in scenario $s$.
	\item[$v_i^{t, s}$] Voltage magnitude at bus $i$ in interval $t$ in scenario $s$.
	\item[$P_{\text{G}, m}^{t, s}$, $Q_{\text{G}, m}^{t, s}$] Aggregated active/reactive power in microgrid $m$ in scenario $s$.
	\item[$P_{\text{DG},m}^{t, s}$, $Q_{\text{DG},m}^{t, s}$] Active/reactive power generation of equivalent dispatchable DG in microgrid $m$ in scenario $s$.
	\item[$P_{\text{d},m}^{t, s}$ and $Q_{\text{d},m}^{t, s}$] Active/reactive load in microgrid $m$ in scenario $s$.
	\item[$E_{\text{DG}, m}^{t, s}$] Energy of equivalent dispatchable DG by the end of interval $t$ in scenario $s$.
\end{IEEEdescription}

\section{Introduction}
\IEEEPARstart{R}{ecent} major blackouts caused by extreme weather events lead to catastrophic consequences for the economy and society \cite{Bie2017, Panteli2017}. 
The impacts of extreme weather events pose unprecedented challenges to power grids and  emphasize the importance of improving system resilience \cite{Wang2016f,  Li2019b, Tabatabaei2019}.
As distribution systems remain vulnerable to natural disasters, it is indispensable to restore the electric service effectively in response to severe power outages, thus achieving more resilient distribution systems \cite{Chen2017}. When severe blackouts occur, a variety of local resources, e.g., microgrids and distributed energy resources (energy storage systems, etc.), can be utilized to restore critical loads in distribution systems. Moreover, the emerging mobile energy storage systems (MESSs) \cite{IEEE2019} can provide temporal-spatial mobility and coordinate with stationary local resources for an integrated distribution system restoration.

Great progress has been made in the utilization of stationary resources for service restoration in distribution systems after major blackouts \cite{Wang2016f, Xu2018}. Microgrids can consolidate and manage a wide range of distributed energy resources to alleviate the hazardous impacts of extended outages \cite{Chen2017}. Reference \cite{Huang2017} proposes a resilience response framework by generator re-dispatch, topology switching, and load shedding. A Markov model is proposed to construct sequential proactive strategies against extreme weather events in \cite{Wang2017d}. In \cite{Amirioun2018}, a microgrid proactive management framework is proposed to coordinate generation reschedule, conservation voltage regulation and demand-side resources. 
Proactive scheduling in multiple energy carrier microgrids is proposed in response to approaching hurricane \cite{Amirioun2018a} and floods \cite{Amirioun2018b}. Reference \cite{Wang2018b} proposes an optimal restoration strategy that coordinates multiple sources at multiple locations to serve critical loads after blackouts. 
These studies illustrate the value of coordination of multiple resources to enhance grid resilience. In addition, with the increasing installation of charging/discharging facilities \cite{KumarNunna2018}, microgrids can provide plug-and-play integration of MESSs for effective service restoration.

MESSs are generally vehicle-mounted container battery energy storage systems with standard interfaces that allow for plug-and-play \cite{IEEE2019}. The importance of the integration of MESS fleets with power system operation has been increasingly recognized in recent researches. For normal operations, MESSs are employed to achieve load shifting \cite{Abdeltawab2017}  and relieve transmission congestion  \cite{Sun2016c}.
In response to extreme events, MESS fleets can be utilized in both pre and post stage. The economic feasibility of MESSs is demonstrated in \cite{Kim2018} by optimizing the investments and relocation of MESSs in case of natural disasters. Reference \cite{Lei2016} proposes a sequential framework for pre-positioning of mobile generators to staging locations and real-time dispatching to distribution systems. In \cite{Gao2017}, the resource allocation of electric buses and transportable batteries is formulated for proactive preparedness for extreme weather events. Dynamic microgrids formation is applied to accommodate mobile and stationary distributed generation and energy resources after disruptions in \cite{Sedzro2018}. Reference \cite{Che2018} presents a microgrid-based critical load restoration by adaptively forming microgrids and positioning mobile emergency resources. Nevertheless, the resource allocation is for one-time dispatching of MESSs in the pre or initial stage of disasters instead of optimizing the temporal-spatial behaviors throughout the restoration process, so the mobility and flexibility of MESS fleets are underutilized. 
Reference \cite{Lei2018} implements resilient routing and scheduling of mobile power sources via a two-stage framework, in which the pre-position and dynamic dispatch are used to coordinate with conventional restoration efforts.
For post-disaster restoration, 
\cite{Lei2018a} proposes a resilient scheme for disaster recovery logistics that involves scheduling of repair crews and mobile power source and network reconfiguration.  A joint scheme is proposed in \cite{Yao2018} to integrate the dynamic scheduling of MESSs, generation re-dispatching and network reconfiguration. However, these researches are either deterministic or do not thoroughly investigate the potential subsequent damage and repair during the restoration process, and more detailed stochastic modeling of MESS in transportation network are still needed.

Furthermore, extreme weather events can destroy not only the distribution systems but also some other interdependent infrastructures \cite{Bie2017, Wang2016f}, 
e.g., transportation networks, which in turn will impact the scheduling of MESSs and impose more challenges to service restoration. 
Few existing studies have considered the electric service restoration in an integrated distribution and transportation system. 
In addition, during the disasters, multiple sources of information can be utilized to improve situational awareness of damage status \cite{Chen2017, Madani2015}, i.e. weather forecast combined with the geographic information systems, distribution system data from smart meters and micro-phasor measurement can provide information on damage and repair to  both the roads in transportation networks  and the branches in distribution systems. Therefore, an integrated restoration strategy is needed to coordinate the mobile and stationary resources for service restoration with dynamically updated system damage information in coupled transportation and distribution systems.


In this context, this paper aims to bridge the gap in the coordination of MESS fleets with microgrids into distribution system restoration and leveraging dynamically updated information during the restoration process. A rolling integrated restoration strategy is proposed to coordinate the dynamic scheduling of MESS, resource dispatching of microgrids and distribution network reconfiguration. In order to take advantages of multiple source data that improves situational awareness during the restoration process, a rolling optimization is adopted to dynamically update system damage status. 
The optimization problem at each interval over the prediction horizon is formulated as a two-stage stochastic mixed-integer linear program (MILP), aiming to minimize the total cost  by co-optimizing the scheduling problem of an MESS fleet, generation dispatching of microgrids and network topology reconfiguration. The contributions of the paper are concluded as follows.

1) A novel integrated restoration strategy is proposed that coordinates the MESS fleet and microgrids to minimize the total system cost. 
The scheduling of MESS fleet is modeled by a stochastic multi-layer time-space network, which reduces the  computational complexity with fewer number of binary variables and constraints and can be utilized for practical transportation networks.

2) The proposed model takes into account both damage and repair to the roads in transportation networks and the branches in distribution systems. The uncertainties in load consumption and the status of the roads and branches are considered to generate scenarios by Monte Carlo simulation method. A rolling optimization framework is adopted to dynamically update system information and the coordinated scheduling over the prediction horizon is formulated as a two-stage stochastic MILP.

The remainder of this paper is organized as follows. Section \ref{sec: Modeling of Mobile Energy Storage Systems} describes the construction of time-space networks and the stochastic  scheduling of MESSs. Section \ref{sec: Rolling optimal integrated restoration strategy} presents the rolling optimization framework for integrated service restoration. Section \ref{sec: Case Studies} conducts case studies on two integrated test systems to verify the effectiveness of the proposed method. Section \ref{sec: Conclusions} summarizes the paper.


\section{Stochastic Modeling of Mobile Energy Storage System}
\label{sec: Modeling of Mobile Energy Storage Systems}

The increasing penetration of MESSs highlights the superiorities over stationary resources in terms of mobility and flexibility. When major blackout happens, MESS can be dispatched among microgrids to transport energy for service restoration. This section formulates the stochastic modeling of MESSs via a time-space network, which has been employed to investigate the vehicle routing and scheduling problem \cite{Yan2014, Lu2018, Liu2018c}. 
In order to take into account the uncertainties in damage and repair to the roads in the transportation networks, a scenario-based stochastic time-space network is proposed to model the temporal-spatial behavior of MESSs over the transportation network, while the charging/discharging schedule is described by battery operation and temporal-spatial constraints. 

\subsection{Construction of Multi-layer Time-space Networks}


A transportation network is modeled as a weighted graph $\mathcal{G}_\text{T}=(\mathcal{N}_\text{T}, \mathcal{E}_\text{T}, \mathcal{W}_\text{T})$, where $\mathcal{N}_\text{T}$ is the nodes set and $\mathcal{E}_\text{T}$ denotes the edges set of roads with the edge distance $\mathfrak{w}  \in \mathcal{W}_\text{T}$.

A set of microgrids $\mathcal{M}$ indexed by $m$ and a set of depots $\mathcal{D}$ are located in the transportation network $\mathcal{G}_\text{T}$. The mappings $\mathfrak{F}_\text{T}:\mathcal{M} \rightarrow \mathcal{N}_\text{T}$ and $\mathfrak{F}_\text{D}:\mathcal{D} \rightarrow \mathcal{N}_\text{T}$ denotes microgrids and depots' locations in the transportation network, respectively. $\Omega$ represents an MESS fleet. An MESS $\omega \in \Omega$ is initially located at a depot $d\in \mathcal{D}$, where it starts and travels among microgrids to provide power supply to power grids, finally it goes back to a depot. 

In order to account for the uncertainties in damage and repair to the roads in transportation networks, a scenario-based stochastic model is adopted. The uncertainty modeling and scenario generation are detailed in Section \ref{sec: Rolling optimal integrated restoration strategy}-B. In a scenario $s$, the shortest path matrix $\mathfrak{P}^s$ is utilized to define the shortest paths for all pairs of microgrids or depots, where the superscript $s$ represents scenario $s$ and the element $ \mathfrak{p}_{ij}^s$ denotes the set of nodes $\mathcal{N}_{ij}^s$  and edges $\mathcal{E}_{ij}^s$  with edge distances $\mathcal{W}_{ij}^s$ in the route, which is calculated by the Dijkstra's algorithm \cite{Cormen2009}. 
It is noted that data-driven approach can be integrated in the future research to find the optimal path considering uncertainties in travel time \cite{Cao2016b, Cao2016c}. 
In this section, a distance matrix $\mathfrak{D}^s$ describes the distances between every two microgrids or depots through the shortest path, the element $\mathfrak{d}_{ij}^s$ is calculated by the sum of the edge distances along the shortest path. A travel time matrix $\mathfrak{T}^{\omega, s}$ with elements $\mathfrak{t}_{ij}^s$ indicates the travel time in the number of intervals considering MESSs' average speed $V_{\text{avg}}^\omega$. 

\begin{gather}
\label{eq: road_matrix_1}
\mathfrak{p}_{ij}^s = (\mathcal{N}_{ij}^s, \mathcal{E}_{ij}^s, \mathcal{W}_{ij}^s),
\forall i,j \in \mathcal{M} \cup \mathcal{D}, s \\
\label{eq: road_matrix_2}
\mathfrak{d}_{ij}^s = \sum_{\mathfrak{w} \in \mathcal{W}_{ij}^s} \mathfrak{w},
\forall s \\
\label{eq: road_matrix_3}
\mathfrak{t}_{ij}^{\omega, s} = \lceil \mathfrak{d}_{ij}^s / V_{\text{avg}}^\omega / \Delta t \rceil, 
\forall \omega \in \Omega, i,j \in \mathcal{M} \cup \mathcal{D}, s
\end{gather}

\noindent where $\lceil~\rceil$  is the ceiling function and $\Delta t$ is the length of time intervals.

A modified multi-layer time-space network is proposed to formulate the vehicle scheduling problem of MESSs over a transportation network. For instance, a transportation network connecting four microgrids and one depot is used for illustration, as shown in \ref{fig: transportation_network}. Time horizon $\mathcal{T}$ is the set of time intervals indexed by $t$. The MESS $\omega$ starts at the depot, and obtains the matrix $\mathfrak{P}^s$, $\mathfrak{D}^s$, $\mathfrak{T}^{\omega, s}$. Then its temporal-spatial behavior will be modeled through a time-space network.

\begin{figure}[!tbp]
	\centering
	
	\includegraphics[width=\columnwidth, clip]{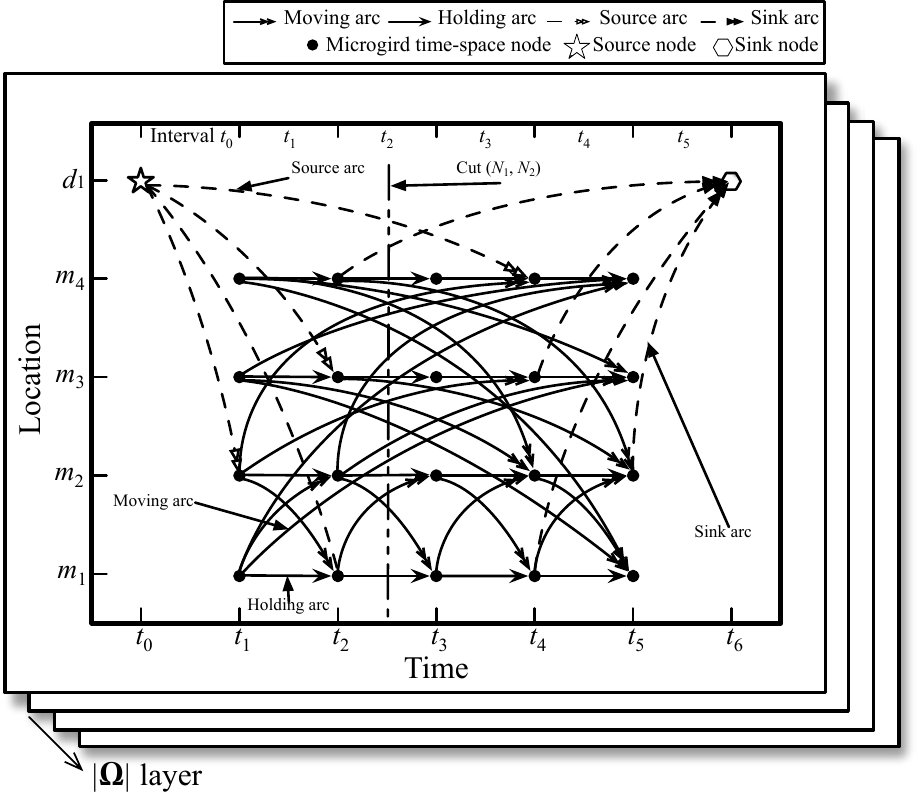}
	\caption{A multi-layer time-space network for modeling temporal-spatial behavior of MESSs over the transportation network.}
	\label{fig: time-space network}
\end{figure}

In the time-space network, as shown in Fig. \ref{fig: time-space network}, the horizontal axis shows the time horizon which is discretized into multiple time intervals, the vertical axis represents the spatial dimension and consists of microgrids and depots. In a scenario $s$, the time-space network is separated into several layers, each one is associated with an MESS $\omega$. That is, $|\Omega|$ layers of the time-space network are assigned to schedule MESSs, where $|\Omega|$ is the total number of MESSs. In a time-space network layer $\mathcal{G}_\text{S}^{\omega,s}=(\mathcal{N}_\text{S}^{\omega,s}, \mathcal{E}_\text{S}^{\omega,s})$ for MESS $\omega$, the set of time-space node $\mathcal{N}_\text{S}^{\omega,s}$ represents microgrids or depots' locations at specific time points. There are three types of time-space nodes as follows. \textit{1) Microgrid time-space nodes $\tilde{\mathcal{N}}_\text{S}^{\omega,s}$}: represent the microgrids at specific time points in the scenario $s$.
\textit{2) Source nodes $\hat{\mathcal{N}}_\text{S}^{\omega,s}$} : represent the MESS $\omega$'s initial depot position in scenario $s$. 
\textit{3) Sink nodes $\check{\mathcal{N}}_\text{S}^{\omega,s}$} : indicate the MESS's final depot positions in scenario $s$.

In addition, the time-space arcs $ \mathcal{E}_\text{S}^{\omega,s}$ connect time-space nodes and describe the feasible movements among locations considering travel time in scenario $s$. Four types of arcs are defined in the time-space network as follows. 
\textit{1) Moving arcs $\tilde{\mathcal{E}}_\text{S}^{\omega,s}$}: a moving arc connects two microgrid time-space nodes and represents a movement in spatial and time dimensions in scenario $s$. For example,  as shown in Fig. \ref{fig: time-space network}, the specified moving arc represents that it is available for an MESS to move from microgrid \#1 at $t_1$ to microgrid \#3 at $t_5$. The movement takes 4 time intervals, which is obtained from the aforementioned travel time matrix $\mathfrak{T}$. All moving arcs that are beyond the time horizon are infeasible and removed from the $\tilde{\mathcal{E}}_\text{S}^{\omega,s}$. In addition, the moving arcs will trigger transportation costs for MESSs.
\textit{2) Holding arcs $\bar{\mathcal{E}}_\text{S}^{\omega,s}$}: a holding arc connects two time-space nodes for the same microgrids in a time interval in scenario $s$. As shown in Fig. \ref{fig: time-space network}, holding arcs indicate that MESSs can stay at microgrids for an interval. Only when an MESS stays on the holding arc can it charge from or discharge to distribution systems.
\textit{3) Source arcs $\hat{\mathcal{E}}_\text{S}^{\omega,s}$}: a source arc connects a source node to a microgrid time-space node in scenario $s$, which implies that MESS $\omega$ is initially located at a depot and moves to a microgrid.
\textit{4) Sink arcs $\check{\mathcal{E}}_\text{S}^{\omega,s}$}: a sink arc connects a microgrid time-space nodes to a sink node in scenario $s$, indicating that MESS $\omega$ return from a microgrid time-space node to a depot at the end of the time horizon. 

\subsection{Impact Analysis of Damage and Repair to Roads on Time-space Network}

\begin{figure}[!tbp]
	\centering
	
	\includegraphics[width=\columnwidth, clip]{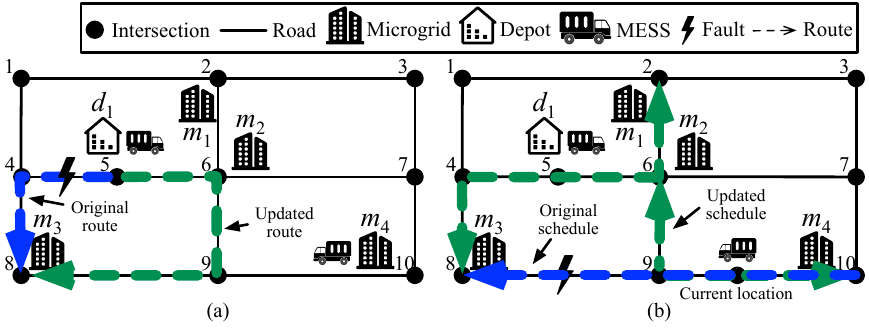}
	\caption{A transportation network connecting microgrids and depots: (a) impacts of damage to roads, (b) rescheduling or rerouting at each interval.}
	\label{fig: transportation_network}
\end{figure}

The damage and repair to the roads in transportation networks have impacts on the matrix $\mathfrak{P}^s$, $\mathfrak{D}^s$, $\mathfrak{T}^{\omega,s}$, leading to the reconstruction of time-space arcs in the time-space network. For example in Fig. \ref{fig: transportation_network}(a), suppose in one of the scenarios the road 4-5 is at fault in the first interval and will be repaired in interval 3. The shortest path from depot \#1 to microgrid \#3 changes due to the road damage, so do the distance and travel time. The original routes and updated routes are denoted as the blue and green dash lines, respectively in Fig \ref{fig: transportation_network}(a). It indicates that an MESS now takes four intervals to travel from depot \#1 to microgrid \#3, compared to two intervals before the road damage. Thus, in the time-space network, the corresponding arcs connecting depot \#1 at time $t_0$ to microgrid \#3 should be modified from $t_2$ to $t_4$. 
Furthermore, after the road 4-5 gets repaired, the travel time is reduced to two intervals and corresponding arcs are modified.

In addition, a rolling optimization framework is utilized to update the scheduling of MESS fleet at each interval. The detailed description of the rolling optimization framework is given in Section \ref{subsec: rolling}. In this formulation, MESSs' current locations are obtained at each interval as the initial condition for the optimization over the new prediction horizon,MESSs can be re-dispatched to any microgrids or depots. 
For instance, an MESS is dispatched from microgrid \#4 to microgrid \#3 via the route 10-9-8, indicated as the blue dash line in Fig .\ref{fig: transportation_network}(b). When it implements the decision for the first interval and gets to a new location, as indicated by the current location. The matrix $\mathfrak{P}^s$, $\mathfrak{D}^s$, $\mathfrak{T}^{\omega, s}$ are updated  to reconstruct the time-space network. Based on the updated initial condition and system information, the MESS does not have to go the microgrid \#3, it can go to other microgrids or even go back to microgrid \#4. 
Moreover, this formulation can also deal with subsequent road damage occurring during the MESS's movement. For example, as shown in Fig. \ref{fig: transportation_network}(b), the damage to road edge 8-9 occurs when the MESS move to microgrid \#3, the designated route is no longer available, the matrix $\mathfrak{P}^s$, $\mathfrak{D}^s$, $\mathfrak{T}^{\omega, s}$ need to be updated considering the road damage to provide the initial condition for the next decision. It is also assumed that if the current road is damaged, the MESS needs to move back to the nearest intersection node with the same road. For example, if the MESS is already on the road edge 8-9 when the damage occurs, it is assumed that the MESS cannot move forward and can only move back to the node 9.
Similarly, the repair can also be taken into account by updating the matrix $\mathfrak{P}^s$, $\mathfrak{D}^s$, $\mathfrak{T}^{\omega, s}$. Therefore, the time-space network is modified based on the updated travel time matrix $\mathfrak{T}^{\omega,s}$.

\subsection{Temporal-spatial Constraints of Mobile Energy Storage Systems}

The temporal-spatial behavior of MESSs over the transportation network $\mathcal{G}_\text{T}$ is transformed into the multi-layer time-space network $\mathcal{G}_\text{S}^s=(\mathcal{N}_\text{S}^s, \mathcal{E}_\text{S}^s)$. The scheduling of vehicles is defined as a sequence of trips by time-space arcs starting from a source node to microgrid time-space nodes and finally returning to a sink node. The formulation is based on arc-wise binary variables $ \zeta _{\hat{n}\check{n}}^{\omega,s}$, which are 1 if the MESS $\omega$ (corresponding to the $\omega^\text{th}$ layer) is on an arc $(\hat{n}, \check{n}) \in \tilde{\mathcal{E}}_\text{S}^{\omega,s}$ in scenario $s$, otherwise set to 0. 

Consider the time-space network in Fig. \ref{fig: time-space network}, a cut $(\mathcal{N}_1^s, \mathcal{N}_2^s)$ of $\mathcal{G}_\text{S}^{\omega,s}=(\mathcal{N}_\text{S}^{\omega,s}, \mathcal{E}_\text{S}^{\omega,s})$ is defined as a partition of $\mathcal{N}_\text{S}^{\omega,s}$ into two disjoint subsets (where $\mathcal{N}_1^s$ represents the nodes on the left side of the cut while $ \mathcal{N}_2^s$ denotes ones on the right). The cut-set of the cut $(\mathcal{N}_1^s, \mathcal{N}_2^s)$ is a set $\{(n_1, n_2) \in  \mathcal{E}_\text{S}^{\omega,s} | n_1\in \mathcal{N}_1^s, n_2\in \mathcal{N}_2^s \}$. Any edge in the cut-set has one endpoint at each side of the cut $(\mathcal{N}_1^s, \mathcal{N}_2^s)$. There is a cut for each time interval in time-space network and an associated cut-set, which is denoted by $\mathbf{C}_\omega^{t,s}$  representing the one for MESS $\omega$ in the interval $t$. For example, a cut for interval $t_2$ is depicted in Fig. \ref{fig: time-space network}. The cut-set $\mathbf{C}_\omega^{t_2}$ involves all the time-space arcs that crossing the cut, indicating all the feasible states of the MESS $\omega$ in the interval and an MESS can only be on exact one arc  in any time interval. The MESS's states can be described by constraint \eqref{eq: tsn_form_1}:

\begin{gather}
\label{eq: tsn_form_1}
\sum_{(\hat{n}, \check{n})\in \mathbf{C}_\omega^{t,s}}{\zeta _{\hat{n}\check{n}}^{\omega, s}} = 1,
\forall \omega, s, t  \\ 
%
\label{eq: tsn_form_2}
\sum_{(\hat{n}, \check{n}) \in \acute{\mathcal{E}}_{\text{S}, n}^{\omega, s}} {\zeta _{\hat{n}\check{n}}^{\omega, s}} = \sum_{ (\hat{n}, \check{n}) \in \grave{\mathcal{E}}_{\text{S}, n}^{\omega, s}} {\zeta _{\hat{n}\check{n}}^{\omega, s}}, \forall n \in \tilde{\mathcal{N}}_\text{S}^{\omega, s}, \omega, s \\
\label{eq: tsn_form_3}
\sum_{(\hat{n}, \check{n}) \in \acute{\mathcal{E}}_{\text{S}, n}^{\omega, s}} {\zeta _{\hat{n}\check{n}}^{\omega, s}} = \zeta _\text{init}^{\omega, s},
\forall n \in \hat{\mathcal{N}}_\text{S}^{\omega, s}, \omega, s
\end{gather}

In addition, for each time-space node $n \in \tilde{\mathcal{N}}_\text{S}^{\omega, s}$ in scenario $s$, the time-space arcs connecting the $n$ can be classified into two groups of in-flows $\grave{\mathcal{E}}_{\text{S}, n}^{\omega, s}$ and out-flows $\acute{\mathcal{E}}_{\text{S}, n}^{\omega, s}$, which represent the arcs entering or leaving the $n$, respectively. For source nodes $n \in \hat{\mathcal{N}}_\text{S}^{\omega, s}$, there are out-flows with the initial location of MESSs. Each time-space node $n$ needs to satisfy the network flow conservation, that is,  an MESS $\omega$ ends trips at time-space node $n$, which serves as the starting point of the subsequent trip.  The $\zeta _\text{init}^{\omega, s}$ indicates the MESS's initial location. The time-space network flow conservation is suggested by constraint \eqref{eq: tsn_form_2}-\eqref{eq: tsn_form_3}.

\subsection{Operation Constraints of Mobile Energy Storage Systems}

When staying on a holding arc, an MESS can exchange power through charging from or discharging to distribution systems. The holding arcs for microgrid $m$  at time interval $t$ in scenario $s$ is denoted by $\bar{\mathcal{E}}_{\text{S}, m}^{\omega, t, s}$, which can be obtained by determining the holding arc that involves time-space nodes for microgrid $m$  in time interval $t$ in scenario $s$. The charging/discharging behaviors are constrained as follows. 

\begin{gather}
 \label{eq: tsn_op_1}
 0\leq P_{\text{ch}, m}^{\omega, t, s} \leq \zeta_{\hat{n}\check{n}}^{\omega, s} \overline{P}_\text{ch}^{\omega},
 \forall (\hat{n}, \check{n}) \in \bar{\mathcal{E}}_{\text{S}, m}^{\omega, t, s}, \omega, m, t, s \\
 \label{eq: tsn_op_2}
 0\leq P_{\text{dch}, m}^{\omega, t, s} \leq \zeta_{\hat{n}\check{n}}^{\omega, s} \overline{P}_\text{dch}^{\omega},
 \forall (\hat{n}, \check{n}) \in \bar{\mathcal{E}}_{\text{S}, m}^{\omega, t, s}, \omega, m, t, s \\
 \label{eq: tsn_op_3}
 0\leq \sum_{m\in \mathcal{M}} P_{\text{ch}, m}^{\omega, t, s} \leq I_\text{ch}^{{\omega}, t, s} \overline{P}_\text{ch}^\omega,
 \forall \omega, t, s \\
 \label{eq: tsn_op_4}
 0\leq \sum_{m\in \mathcal{M}} P_{\text{dch}, m}^{\omega, t, s} \leq I_\text{dch}^{{\omega}, t, s} \overline{P}_\text{dch}^\omega,
 \forall \omega, t, s \\
 \label{eq: tsn_op_5}
 I_\text{dch}^{\omega, t, s} + I_\text{dch}^{\omega, t, s} \leq \sum_{(\hat{n}, \check{n})\in \bar{\mathcal{E}}_{\text{S}, m}^{\omega, t, s}} \zeta _{\hat{n}\check{n}}^{\omega, s}, 
 \forall m, \omega, t, s \\
\label{eq: tsn_op_6}
\begin{align} 
E_\omega^{t+1, s} = E_\omega^{t, s} - \Delta t \Big(\frac{\sum_{m\in \mathcal{M}} P_{\text{dch}, m}^{\omega, t, s} }{\eta_\text{dch}^\omega} - \eta_\text{ch}^\omega \sum_{m\in \mathcal{M}} P_{\text{ch},m}^{\omega, t+1, s}\Big), \notag \\
\forall \omega, t\in \mathcal{T} \setminus\{T\}, s 
\end{align} \\
\label{eq: tsn_op_7}
E_\text{c}^\omega \uline{SOC}^\omega \leq E_{\omega}^{t, s} \leq E_\text{c}^\omega \overline{SOC}^\omega, \forall \omega, t, s
\end{gather}

\noindent Constraints \eqref{eq: tsn_op_1}-\eqref{eq: tsn_op_2} state the relation between charging/discharging power and temporal-spatial behaviors. Equations \eqref{eq: tsn_op_3}-\eqref{eq: tsn_op_4} express the charging/discharging power associated with battery status, which is also constrained by temporal-spatial behaviors in \eqref{eq: tsn_op_5}. Constraint \eqref{eq: tsn_op_6} calculates the energy and constraint \eqref{eq: tsn_op_7} sets the upper and lower range.

\subsection{Complexity Analysis of Time-space Network}

References \cite{Sun2017, Yao2018} add virtual nodes to form time-space arcs that span more than one interval. The time-space network can be extended to an extremely large model, leading to high computational complexity when applied to practical transportation networks. For instance, it takes four time spans from microgrid \#3 to microgrid \#4, three virtual nodes need to be added to represent the time-space arcs from microgrid \#3 at $t$ to microgrid \#4 at $t+4$. Then the arcs between virtual nodes and original nodes are added. This would significantly increase the numbers of binary variables and constraints.

This paper proposes a time-space network without virtual nodes for the arcs that span more than one interval, resulting in fewer numbers of binary variables and constraints. The number of virtual nodes is indicated by $\Delta \mathcal{N}_\text{V}^\omega$, which is determined by the travel time matrix $\mathfrak{T}^{\omega, s}$, the comparison of numbers of variables and constraints are presented in Table \ref{table: number_comparison}, where the $P$ is the permutation. Take the case in Fig. \ref{fig: time-space network} and Fig. \ref{fig: transportation_network} for instance, suppose the time horizon is 6 hours and the length of interval is 1 hour. For this scenario, the proposed time-space network reduces the numbers of binary variables and constraints by 58.01\% and 67.53\%, respectively.

\begin{table}[!tbp]
	\centering
	\caption{Comparison of numbers of variables and constraints} \label{table: number_comparison}
	\resizebox{1.03\columnwidth}{!}{
	\begin{tabular}{lll}
		\hline
		& Ref. \cite{Sun2017, Yao2018} & Proposed model \\ \hline
		\begin{tabular}[c]{@{}l@{}}\# of Binary \\ Variables\end{tabular} 
			& \begin{tabular}[c]{@{}l@{}}$\sum_{s \in \mathcal{S}, \omega \in \Omega} |\mathcal{E}_\text{S}^{\omega, s}|+$ \\ $\sum_{s \in \mathcal{S}, \omega \in \Omega} [(|\Delta \mathcal{N}_\text{V}^{\omega, s}|+1)$  \\
				$ \times 2|\mathcal{T}|] $  \\ $   -P_2^{|\mathcal{M} \cup \mathcal{D}|} \times |\mathcal{S}|$ \end{tabular} 
			&    $\sum_{s \in \mathcal{S}, \omega \in \Omega} |\mathcal{E}_\text{S}^{\omega, s}|$         \\ \hline
		\begin{tabular}[c]{@{}l@{}}\# of \\ Constraints\end{tabular}      
			& \begin{tabular}[c]{@{}l@{}}$\sum_{s \in \mathcal{S}, \omega \in \Omega} |\mathcal{N}_\text{S}^{\omega, s}|+$ \\ $\sum_{s \in \mathcal{S}, \omega \in \Omega} (|\Delta \mathcal{N}_\text{V}^{\omega, s}| \times |\mathcal{T}|) $\\ \end{tabular}    
			&     $\sum_{s \in \mathcal{S}, \omega \in \Omega} |\mathcal{N}_\text{S}^{\omega, s}|$              \\ \hline
	\end{tabular}
	}
\end{table}

\section{Rolling Integrated Service Restoration}
\label{sec: Rolling optimal integrated restoration strategy}

This section proposes a rolling optimal  service restoration that coordinates the scheduling of MESS fleets, resource dispatching of microgrids and distribution network reconfiguration. The objective is to minimize the total cost, considering the customer interruption cost, microgrid generation cost, and MESS transportation cost and battery maintenance cost. Partial load curtailment is used as in \cite{Ding2017a}. Customer interruption cost is adopted to differentiate between critical and non-critical loads \cite{Li2014}. In order to take into account the subsequent damage during the restoration process in both distribution systems and transportation networks, a rolling optimization is adopted for dynamic updating of system damage status. The detailed problem statement, rolling optimization framework, and mathematical formulation are described as follows.

\subsection{Problem Statement}

A conceptual resilience curve associated with an event in \cite{Panteli2015a} is adopted for better illustration, as shown in Fig. 1. $R$ refers to an index of system resilience level. The system states involve: pre-disturbance resilient state ($t_0$, $t_\text{e}$), event progress ($t_\text{e}$, $t_\textit{pe}$), post-event degraded state ($t_\textit{pe}$, $t_\text{r}$), restorative state ($t_\text{r}$, $t_\text{pr}$), post-restoration state ($t_\text{pr}$, $t_\text{ir}$) and infrastructure recovery ($t_\text{ir}$, $t_\text{pir}$). 
%
%
It is assumed in this manuscript that extreme events cause the complete outages of transmission grids and the distribution systems can no longer be supplied by transmission grids. Under this circumstance, microgrids can be utilized to coordinate multiple stationary and mobile resources for service restoration at distribution level. It is noted that each load is powered by only one microgrid \cite{Xu2018, Gao2017, Lei2016}, there is no loop or overlap region. However, the model can be further extended to consider control strategies in network reconfiguration \cite{Bie2017}.

\begin{figure}[!tbp]
	\centering
	\includegraphics[width=\columnwidth, clip]{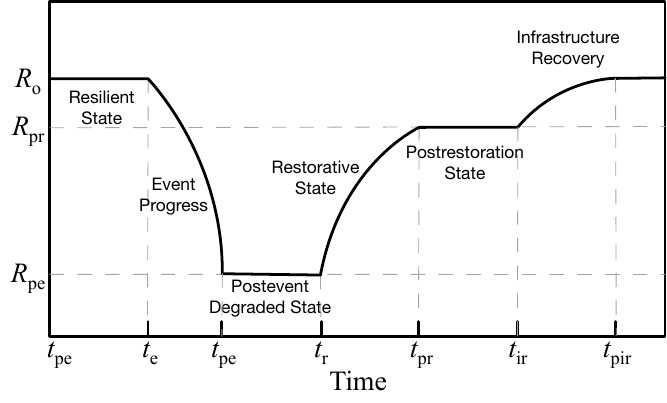}
	\caption{Conceptual resilience curve associated with an event \cite{Panteli2015a}.}
	\label{fig: resilience_curve}
\end{figure}

Previous studies generally assume that all the faults have been identified at $t_\text{r}$ \cite{Yao2018, Xu2018} and the restoration starts with accurate damage information. However, the extreme events may involve multiple stages, i.e., after the major strikes, the extreme events can still last and cause subsequent damage. Therefore, the proposed integrated restoration strategy is implemented right after the major strikes of an extreme event until the main grid is restored, i.e., from $t_\text{pe}$ to $t_\text{ir}$, to enhance the system resilience level. In this case, restoration gets started with incomplete damage information. A two-stage stochastic optimization is adopted to account for uncertainties and a rolling optimization framework is used to dynamically update system information over the prediction horizon at each interval.

\subsection{Uncertainty Modeling and Scenario Generation}
\label{subsec: rolling}

%
%
A scenario-based method is used to model uncertainties, by generating a large number of scenarios and doing scenario reduction to ensure the computational tractability.

There are several uncertainties considered in this paper, including forecasting errors in load, the status of the roads in transportation networks and the branches in distribution systems. A normal distribution is used to represent the forecasting error of load, in which, the mean value of the normal distribution is the predicted load and the standard deviation is set to be 2\% of the predicted load \cite{Wang2014a}. A two-state continuous-time Markov model is applied to represent the availability of roads and branches \cite{Wu2007}, while available and unavailable hours are subject to exponential distributions with mean uptime and mean downtime values \cite{Sun2017}. The repair efforts are also taken into account during the optimization time horizon. Then Monte-Carlo simulation is utilized to generate scenarios.
In order to reduce the computation efforts, a scenario reduction is implemented to reduce the number of scenarios while maintaining a good approximation of the system uncertainty. This paper adopts a simultaneous backward reduction method \cite{Wu2007} for scenario reduction.

\subsection{Rolling Optimization Framework}
During the restoration process, it is challenging to obtain the accurate load forecast and damage status at the very beginning and to solve the optimization problem only once to acquire the acceptable solution. In addition, the load forecast and damage status are dynamically updated from multiple sources of information, with the short-term forecast being more accurate. Therefore, in order to deal with the inaccuracy of forecast over the long horizon and leverage dynamically updated forecasts, a rolling optimization framework is adopted to solve the problem recursively in a finite-moving-horizon of intervals \cite{OConnell2014a}.
Specifically, the entire time horizon $\mathcal{T}_\text{H}$ is discretized into equal time intervals by $\Delta t$ and the optimization problem is formulated and solved at each time interval over the prediction horizon $\mathcal{T}_\text{P}$, but only the decisions in the first interval  are implemented. Then the prediction horizon is shifted forward and the calculation is repeated for the new prediction horizon until the end of the entire time horizon $\mathcal{T}_\text{H}$, based on the updated system information and initial conditions \cite{Li2018}. The final decisions are the sequence of the decisions in the first interval of each prediction horizon.

%


The forecasts for uncertainties in each scenario $s$ are denoted with superscript $s$, while the realization of uncertain parameters is indicated with superscript $\xi$. 
The proposed rolling optimization framework adopts a two-stage stochastic programming approach, in which the decision variables are divided into two different groups: \textit{1)} First-stage variables  involve the decision variables $\zeta _{\hat{n}\check{n}}^{\omega, s}$, $\alpha_{ij}$ in the first interval of prediction horizon $\mathcal{T}_\text{P}$, and are determined before the realization of uncertain parameters and scenario independent by nonanticipativity constraint, as described in the next subsection. 2) The second-stage decisions are scenario-dependent and can be adjusted once uncertain parameters reveal, consisting of decision variables $\zeta _{\hat{n}\check{n}}^{\omega}$ in the rest intervals of $\mathcal{T}_\text{P}$ and $P_{\text{ch}, m}^{\omega, t, s}$, $P_{\text{dch}, m}^{\omega, t, s}$, $P_{\text{DG}, m}^{t, s}$, $Q_{\text{DG}, m}^{t, s}$ in all intervals of $\mathcal{T}_\text{P}$. The solution of the two-stage stochastic program is a single first-stage policy and a collection of second-stage recourse decisions defining recourse actions in response to each scenario $s$. 

The first-stage decisions $\zeta _{\hat{n}\check{n}}^{\omega}$, $\alpha_{ij}$ are implemented at the beginning of the interval $t$. 
At the end of the interval $t$, the actual decisions $P_{\text{ch}, m}^{\omega, t}$, $P_{\text{dch}, m}^{\omega, t}$, $P_{\text{DG}, m}^{t}$, $Q_{\text{DG}, m}^{t}$ can be obtained by a deterministic re-optimization based on the implementation of $\zeta _{\hat{n}\check{n}}^{\omega}$, $\alpha_{ij}$ and realization of uncertainties of $P_{\text{D}, i}^{t, \xi}$, $Q_{\text{D}, i}^{t, \xi}$, $\ddot{\mathcal{E}}_\text{D}^{\xi}$, $\ddot{\mathcal{E}}_\text{T}^{ \xi}$. 
Then the rolling optimization proceeds to the next time interval $t+1$, the initial condition will be updated by the realization of uncertainties of $P_{\text{D}, i}^{t, \xi}$, $Q_{\text{D}, i}^{t, \xi}$, $\ddot{\mathcal{E}}_\text{D}^{\xi}$, $\ddot{\mathcal{E}}_\text{T}^{ \xi}$ and actual decisions of $\zeta _{\hat{n}\check{n}}^{\omega}$, $\alpha_{ij}$, $P_{\text{ch}, m}^{\omega, t}$, $P_{\text{dch}, m}^{\omega, t}$, $P_{\text{DG}, m}^{t}$, $Q_{\text{DG}, m}^{t}$ in the interval $t$.

%
%




\subsection{Mathematical Formulation}
A distribution network is modeled as a graph $\mathcal{G}_\text{D}=(\mathcal{N}_\text{D}, \mathcal{E}_\text{D})$ \cite{Xu2018}, where $\mathcal{N}_\text{D}$ is the set of distribution system buses, indexed by $i$ and $\mathcal{E}_\text{D}$ is the set of distribution system branches, indexed by $(i, j)$.  The mapping $\mathfrak{G}_\text{D}:\mathcal{M} \rightarrow \mathcal{N}_\text{D}$ indicates the microgrid locations in the distribution network.


The network reconfiguration is formulated by a fictitious network model \cite{Ding2017a, Bie2017, Liu2018b} to describe the spanning forest constraints.
A linearized DistFlow model \cite{Lei2016, Ding2017a} is employed for power flow analysis. 
It is noted that the linearized three-phase power flow \cite{Wang2018b} can be included for further extension of this model to unbalanced three-phase conditions. The mathematical formulations are described as follows.

%
%
%
%

\begin{gather}
\label{eq: ds_topo_1}
\sum_{i \in \delta^s (j)} \mathfrak{f}_{ji}^{s}-\sum_{i \in \pi^s (j)} \mathfrak{f}_{ij}^{ s} = -1, 
\forall j\in \mathcal{N}_\text{D} \setminus \mathcal{G}_\text{D}(\mathcal{M}), s \\ 
\label{eq: ds_topo_2}
\sum_{i \in \delta^s (j)} \mathfrak{f}_{j i}^{s}-\sum_{i \in \pi^s (j)} \mathfrak{f}_{i j}^{s}=\mathfrak{h}_{j}^s, 
\forall	j \in \mathcal{G}_\text{D}(\mathcal{M}), s \\ 
\label{eq: ds_topo_3}
{-M \alpha_{i j}^s \leq \mathfrak{f}_{i j}^s \leq M \alpha_{i j}^s},
\forall (i, j), s \\ 
\label{eq: ds_topo_4}
{-M(2-\alpha_{i j}^s) \leq \mathfrak{f}_{i j}^s \leq M(2-\alpha_{i j}^s)}, 
\forall (i,j), s \\ 
\label{eq: ds_topo_5}
{\mathfrak{h}_{j}^s \geq 1}, \forall j \in  \mathcal{G}_\text{D}(\mathcal{M}), s\\
\label{eq: ds_topo_6}
\sum_{(i,j)\in \mathcal{E}_\text{D}} \alpha_{ij}^s = \mathcal{|N_\text{D}|}-\mathcal{|M|}, 
\forall s \\
\label{eq: ds_topo_7}
\alpha_{ij}^s = 0, \forall (i,j) \in \ddot{\mathcal{E}}_\text{D}^s, s\\
\label{eq: op_ds_1}
P_{\text{G}, i}^{t, s} - P_{\text{r}, i}^{t, s} = \sum_{(i,j)\in \mathcal{E}_\text{D}} P_{ij}^{t, s} - \sum_{(k,i)\in \mathcal{E}_\text{D}} P_{ki}^{t, s}, \forall i, t, s \\
 \label{eq: op_ds_2}
Q_{\text{G}, i}^{t, s} - Q_{\text{r}, i}^{t,s } = \sum_{(i,j)\in \mathcal{E}_\text{D}} Q_{ij}^{t, s} - \sum_{(k,i)\in \mathcal{E}_\text{D}} Q_{ki}^{t, s}, \forall i, t, s \\
\label{eq: op_ds_3}
v_j^{t,s} - v_i^{t, s} \leq M(1-\alpha_{ij}^s) + \frac{r_{ij}P_{ij}^{t, s} + x_{ij}Q_{ij}^{t,s}}{v_0}, \forall (i,j), t, s \\
\label{eq: op_ds_4}
v_j^{t, s} - v_i^{t, s} \geq -M(1-\alpha_{ij}^s) + \frac{r_{ij}P_{ij}^{t, s} + x_{ij}Q_{ij}^{t, s}}{v_0}, \forall (i,j), t, s \\
\label{eq: op_ds_5}
-\alpha_{ij}^s \overline{S}_{ij} \leq P_{ij}^{t, s} \leq \alpha_{ij}^s \overline{S}_{ij}, \forall (i,j), t, s \\
\label{eq: op_ds_6}
-\alpha_{ij}^s \overline{S}_{ij} \leq Q_{ij}^{t, s} \leq \alpha_{ij}^s \overline{S}_{ij}, \forall (i,j), t, s \\
\label{eq: op_ds_7}
-\sqrt{2}\alpha_{ij}^s \overline{S}_{ij} \leq P_{ij}^{t, s} + Q_{ij}^{t, s} \leq \sqrt{2}\alpha_{ij}^s \overline{S}_{ij}, 
\forall (i,j), t, s \\
\label{eq: op_ds_8}
-\sqrt{2}\alpha_{ij}^s \overline{S}_{ij} \leq P_{ij}^{t, s} - Q_{ij}^{t, s} \leq \sqrt{2}\alpha_{ij}^s \overline{S}_{ij}, 
\forall (i,j), t, s \\
\label{eq: op_ds_9}
\uline{v}_i \leq v_i^{t, s} \leq \overline{v}_i,
\forall i\in \mathcal{N}_\text{D} \setminus \mathcal{G}_\text{D}(\mathcal{M}), t, s \\
\label{eq: op_ds_10}
v_i^{t,s} = v_0, \forall i\in \mathcal{G}_\text{D}(\mathcal{M}), t, s \\
\label{eq: op_ds_11}
0 \leq P_{\text{r},i}^{t, s} \leq P_{\text{D}, i}^{t, s}, \forall i, t, s \\
\label{eq: op_ds_12}
Q_{\text{r},i}^{t, s} = P_{\text{r},i}^{t, s} \tan(\cos^{-1} \varphi_i), \forall i, t, s \\
\label{eq: op_mg_1}
P_{\text{G}, m}^{t, s} = P_{\text{DG}, m}^{t, s} - \sum_{\omega\in \Omega} (P_{\text{ch}, m}^{\omega, t, s} + P_{\text{dch}, m}^{\omega, t, s}) - P_{\text{d},m}^{t, s},  \forall m, t, s \\
\label{eq: op_mg_2}
Q_{\text{G}, m}^{t, s} = Q_{\text{DG}, m}^{t, s} - Q_{\text{d},m}^{t, s}, \forall m, t, s \\
\label{eq: op_mg_3}
0\leq P_{\text{DG}, m}^{t, s} \leq \overline{P}_{\text{DG}, m}, \forall m, t, s \\
\label{eq: op_mg_4}
-\overline{Q}_{\text{DG}, m} \leq Q_{\text{DG}, m}^{t, s} \leq \overline{Q}_{\text{DG}, m}, \forall m, t, s \\
\label{eq: op_mg_5}
E_{\text{DG},m}^{t+1, s} = E_{\text{DG},m}^{t, s} - P_{\text{DG},m}^{t+1, s}\Delta t,
\forall m, t\in \mathcal{T}_\text{P} \setminus \{T\}, s\\
\label{eq: op_mg_6}
\uline{E}_{\text{DG},m} \leq E_{\text{DG},m}^{t, s} \leq \overline{E}_{\text{DG},m},
\forall m, t, s
\end{gather}




\noindent where $|\mathcal{N}_\text{D}|, |\mathcal{M}|$ represent the cardinality of the sets, $\delta^s (j)$ and $\pi^s (j)$ are the set of children nodes and parent nodes of bus $j$ in scenario $s$, respectively. $\mathfrak{f}_{ij}$ is the power transferred on the line $(i, j)$ in the fictitious network; $\mathfrak{h}_j$ is the power supplied by the “source” buses in the fictitious network. $M$ is a large number. Since the fictitious network has the same topology structure as the original power network, they have the same connectivity. Equations \eqref{eq: ds_topo_1}-\eqref{eq: ds_topo_5} indicates that the satisfaction of energy balance at each bus in the fictitious network implies at least one path exists between the “source” bus and all other buses, so that the sub-graph must be connected. Equation \eqref{eq: ds_topo_6} guarantees the necessary condition for radiality. Equation \eqref{eq: ds_topo_7} represents the damaged branch status.
Constraints \eqref{eq: op_ds_1} and \eqref{eq: op_ds_2} describe the active and reactive power balance at bus $i$. Constraints \eqref{eq: op_ds_3} and \eqref{eq: op_ds_4} indicate the branch voltage drop by the big-M method \cite{Dorostkar-Ghamsari2016}. Equations \eqref{eq: op_ds_5}-\eqref{eq: op_ds_8} provide a linearized approximation regarding branch capacity. Equation \eqref{eq: op_ds_9} suggests the range of voltage magnitude. Equation \eqref{eq: op_ds_10} sets the voltage of microgrid buses to $v_0$. Equations \eqref{eq: op_ds_11} and \eqref{eq: op_ds_12} constrain the load restoration and power factor. Constraints \eqref{eq: op_mg_1} and \eqref{eq: op_mg_2} express the aggregated real and reactive power considering the charging or discharging of MESS fleets. Equations \eqref{eq: op_mg_3} and \eqref{eq: op_mg_4} depict the  power capacity constraints of equivalent dispatchable DG. Equation \eqref{eq: op_mg_5} calculates the energy in each microgrid. Equation \eqref{eq: op_mg_6} presents the range of energy.

Considering that the first stage variables are scenario independent, the nonanticipativity constraints are enforced to ensure that all the realizations of the first-stage decision variables are equal to each other \cite{Shapiro2014}. The nonanticipativity constraints are described as follows.

\begin{gather}
\label{eq: nonanticipativity_1}
\zeta _{\hat{n}\check{n}}^{\omega, s} = \sum_{\mathfrak{s} \in \mathcal{S}} \gamma_\mathfrak{s} \zeta _{\hat{n}\check{n}}^{\omega, \mathfrak{s}}, 
\forall (\hat{n}, \check{n}) \in \mathbf{C}_\omega^{0,s}, \omega, s \\
\label{eq: nonanticipativity_2}
\alpha_{ij}^s = \sum_{\mathfrak{s} \in \mathcal{S}} \gamma_\mathfrak{s} \alpha_{ij}^\mathfrak{s},
\forall (i,j), s 
\end{gather}

\noindent where $\mathbf{C}_\omega^{0,s}$ denotes the cut set of time-space arcs in the first interval of the prediction horizon $\mathcal{T}_\text{P}$ in the scenario $s$.






The objective function is formulated as follows to minimize the total cost.

\begin{align} \label{eq: objective function} 
\min~\sum_{t\in \mathcal{T}_\text{P}} \gamma_s \sum_{s \in \mathcal{S}}   \bigg[
\Big[\sum_{i\in \mathcal{N}} W_i (P_{\text{D},i}^{t, s} - P_{\text{r},i}^{t, s}) + \sum_{m\in \mathcal{M}} C_{\text{gen},m}P_{\text{DG},m}^{t, s} \notag \\ 
+ \sum_{\omega\in \Omega} C_{\text{bat},\omega} \sum_{m\in {\mathcal{M}}}  (P_{\text{ch},m}^{\omega, t, s} + P_{\text{dch}, m}^{\omega, t, s}) \Big] \Delta t \notag \\ 
+   \sum_{\omega\in \Omega} C_{\text{tran},\omega} \sum_{(\hat{n}, \check{n}) \in 
	\mathcal{E}_\text{S}^{\omega, s} } \zeta_{\hat{n}\check{n}}^{\omega, s} \bigg]
\end{align}

\noindent where the term $\sum_{t\in \mathcal{T}} \sum_{s \in \mathcal{S}} \sum_{i\in \mathcal{N}} W_i (P_{\text{D},i}^{t, s} - P_{\text{r},i}^{t, s}) \Delta T$ is the customer interruption cost. The term $\sum_{t\in \mathcal{T}} \sum_{s \in \mathcal{S}} \sum_{m\in \mathcal{M}} C_{\text{gen},m} P_{\text{DG},m}^{t, s}\Delta T$ shows the microgrids generation cost. The third term $\sum_{t\in \mathcal{T}} \sum_{s \in \mathcal{S}}  \sum_{\omega\in \Omega} C_{\text{bat},\omega} \sum_{m\in \mathcal{M}}  (P_{\text{ch},m}^{\omega, t, s} + P_{\text{dch}, m}^{\omega, t, s}) \Delta T$ calculates the MESS battery maintenance cost. The last term $\sum_{s \in \mathcal{S}} \sum_{\omega\in \Omega} C_{\text{tran},\omega} \sum_{(\hat{n}, \check{n}) \in 	\mathcal{E}_\text{S}^{\omega, s}} \zeta_{\hat{n}\check{n}}^{\omega}$ denotes the transportation cost.

The framework of the integrated restoration strategy is illustrated in Algorithm~\ref{algo: dynamic_restoration}.

\begin{figure}[!tbp]
	\removelatexerror
	\begin{algorithm}[H] 
		\caption{Framework for integrated service restoration}  
		\label{algo: dynamic_restoration}
		
		\nonl $\vartriangleright$ \textit{\textbf{Initialization}}:
		
		Input the distribution system $\mathcal{G}_\text{D}=(\mathcal{N}_\text{D},  \mathcal{E}_\text{D})$, transportation network $\mathcal{G}_\text{T}=(\mathcal{N}_\text{T}, \mathcal{E}_\text{T})$, microgrids $m \in \mathcal{M}$ and depots $d \in \mathcal{D}$, mobile energy storage systems $\omega \in \Omega$\;
		
		Generate scenarios and do scenario reduction to obtain $\mathcal{S}$\;
		
		\nonl $\vartriangleright$ \textit{\textbf{Rolling optimization}}:
		
		\For{each interval $t \in \text{time horizon }\mathcal{T}_\text{H}$} 
		{	Move forward the prediction horizon $\mathcal{T}_\text{P}$ \;
			Update initial condition by the realization of uncertainties of $P_{\text{D}, i}^{t-1, \xi}$, $Q_{\text{D}, i}^{t-1, \xi}$, $\ddot{\mathcal{E}}_\text{D}^{\xi}$, $\ddot{\mathcal{E}}_\text{T}^{ \xi}$ and actual decisions of $\zeta _{\hat{n}\check{n}}^{\omega}$, $\alpha_{ij}$, $P_{\text{ch}, m}^{\omega, t-1}$, $P_{\text{dch}, m}^{\omega, t-1}$, $P_{\text{DG}, m}^{t-1}$, $Q_{\text{DG}, m}^{t-1}$\ in the previous interval;

			\For{each scenario $s \in \mathcal{S}$}
			{
				Update the set of damaged branches $\ddot{\mathcal{E}}_\text{D}^s$ and and the set of damaged roads $\ddot{\mathcal{E}}_\text{T}^s$\;
				
				\For{each MESS $\omega \in \Omega$}
				{
					
					Update MESS $\omega$'s current location\;
					
					Compute path matrix $\mathfrak{P}^{\omega,s}$, distance matrix $\mathfrak{D}^{\omega,s}$ and time matrix $\mathfrak{T}^{\omega,s}$ via \eqref{eq: road_matrix_1} - \eqref{eq: road_matrix_3}, 
					
					Construct time-space network $\mathcal{G}_\text{S}^{\omega, s}=(\mathcal{N}_\text{S}^{\omega, s}, \mathcal{E}_\text{S}^{\omega, s})$\;
				}

			}
			
			The optimization problem over the prediction horizon $\mathcal{T}_\text{P}$ is formulated as a two-stage stochastic MILP and solved by CPLEX:
			\vspace{8 pt}
			\nonl \begin{align}
			\min~&\eqref{eq: objective function} \notag \\
			\text{s.t.}~&\eqref{eq: tsn_form_1} -\eqref{eq: nonanticipativity_2} \notag
			\end{align}
			
			\nl Implement optimal solution $\zeta _{\hat{n}\check{n}}^{\omega^\star}$, $\alpha_{ij}^\star$ in the first interval of $\mathcal{T}_\text{P}$\;
			
			
			\nl Actual decisions for $P_{\text{ch}, m}^{\omega, t}$, $P_{\text{dch}, m}^{\omega, t}$, $P_{\text{DG}, m}^{t}$, $Q_{\text{DG}, m}^{t}$  are obtained by a deterministic re-optimization based on the implementation of $\zeta _{\hat{n}\check{n}}^{\omega}$, $\alpha_{ij}$ and realization of uncertainties of $P_{\text{D}, i}^{t, \xi}$, $Q_{\text{D}, i}^{t, \xi}$, $\ddot{\mathcal{E}}_\text{D}^{\xi}$, $\ddot{\mathcal{E}}_\text{T}^{ \xi}$ in interval $t$\;
			
		}
		
		\textbf{return} the solution $\zeta _{\hat{n}\check{n}}^{\omega^\star}$, $P_{\text{ch}, m}^{\omega, t^\star}$, $P_{\text{dch}, m}^{\omega, t^\star}$, $\alpha_{ij}^\star$, $P_{\text{DG}, m}^{t^\star}$, $Q_{\text{DG}, m}^{t^\star}$ for the entire time horizon $\mathcal{T}_\text{H}$.
		
	\end{algorithm}  
\end{figure}

\section{Case Studies}
\label{sec: Case Studies}
The simulations are studied on two integrated distribution and transportation systems to verify the effectiveness of the proposed service restoration strategy, one is with a Sioux Falls transportation network \cite{LeBlanc1975} and four 33-bus distribution systems \cite{Baran1989a}, the other is based on the Singapore transportation network and six 33-bus distribution systems. The proposed model is implemented using Python 3.6 and the two-stage stochastic MILP at each interval over the prediction horizon is solved by CPLEX 12.8.0 \cite{IBM2018a}.

\subsection{Case \Romannum{1}: Sioux Falls Transportation Networks with four 33-bus Distribution Systems}

\vspace{0pt}
\noindent \textit{1) Test Systems}
\vspace{6pt}

\begin{figure}[!tbp]
	\centering
	\includegraphics[width=\columnwidth, clip]{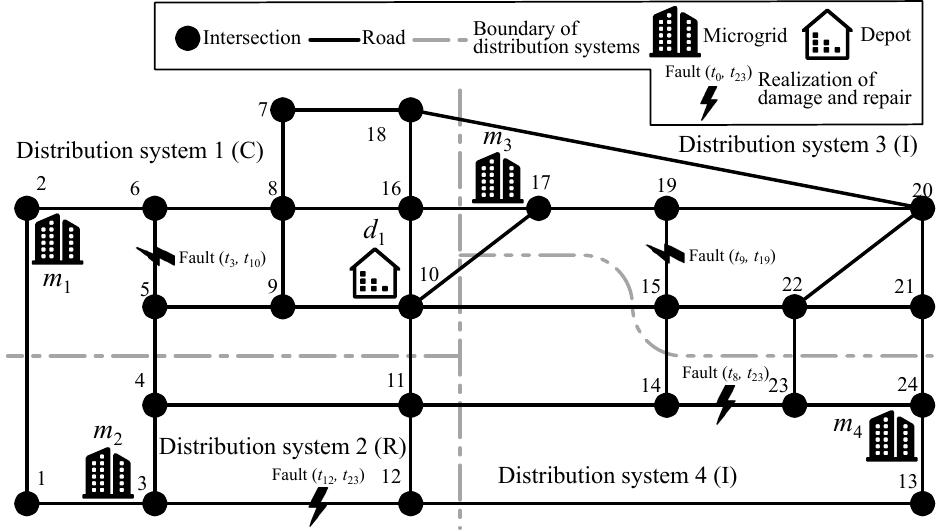}
	\caption{An integrated test system with a Sioux Falls transportation network connecting four distribution systems.}
	\label{fig: integrated_test_system}
\end{figure}

\begin{figure}[!tbp]
	\centering
	\includegraphics[width=\columnwidth, clip]{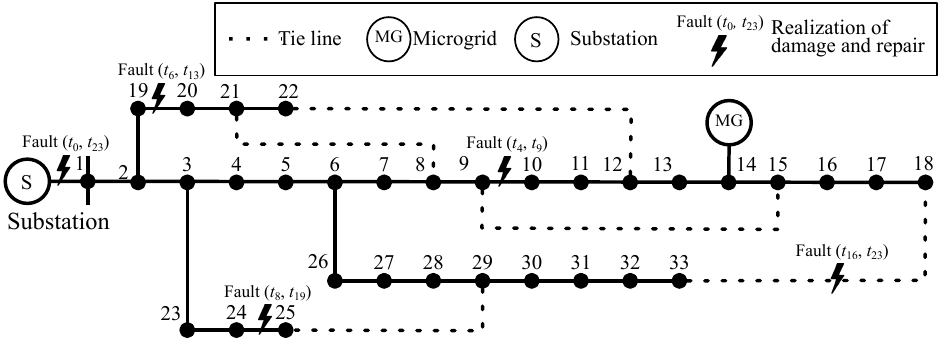}
	\caption{Distribution system: a modified 33-bus test system.}
	\label{fig: 33_bus_system}
\end{figure}

\begin{figure}[!tbp]
	\centering
	\includegraphics[width=\columnwidth, clip]{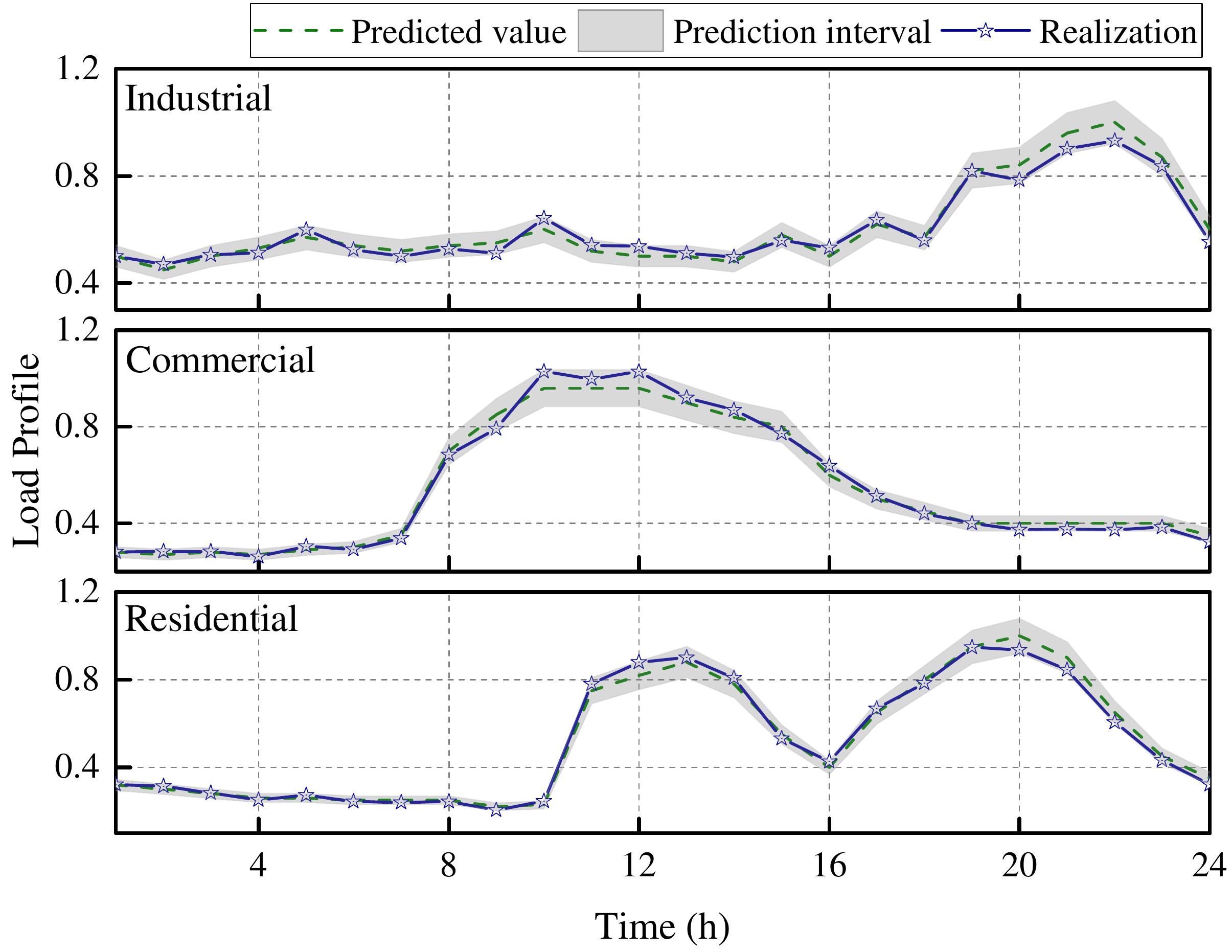}
	\caption{Load profile.}
	\label{fig: load_profile}
\end{figure}

As the scheduling of MESS fleets involves multiple distribution systems, an integrated test system with four 33-bus distribution systems connected by the Sioux Falls transportation network is proposed, as shown in Fig. \ref{fig: integrated_test_system}. The distance of each edge in the transportation network is double. 
As the interfaces between the transportation network and distribution systems are microgrids, other buses and branches of the distribution systems can not directly be linked to the transportation network. Thus, only microgrids' locations are explicitly indicated in Fig. \ref{fig: integrated_test_system} and the boundary of distribution systems are drawn for illustration.
Each distribution system has identical topologies (shown in Fig. \ref{fig: 33_bus_system}) but different categories of loads such as residential (R), commercial (C) and industrial (I). 
Monte Carlo simulation is used to generate 2000 random scenarios, which are reduced to 10 scenarios as described in Section \ref{sec: Rolling optimal integrated restoration strategy}-B. The realization of damage and repair to the roads in the transportation network and the branches in distribution systems are depicted in Fig. \ref{fig: integrated_test_system} and \ref{fig: 33_bus_system}, respectively.
Fig. \ref{fig: load_profile} shows predicted value of industrial, commercial and residential loads, as well as prediction intervals and actual realization. 

The length of entire time horizon $\mathcal{T}_\text{H}$ is set to 24-h while the length of prediction horizon $\mathcal{T}_\text{P}$ is 12-h. 
Unit interruption costs are adopted to distinguish critical and non-critical loads, with $\$10/\text{kWh}$ and $\$2/\text{kWh}$, respectively. Critical loads are randomly selected. 
The microgrid unit generation cost is \$0.5/kWh. The unit battery maintenance cost is $\$0.2/\text{kWh}$. The unit transportation cost is $\$80/\text{h}$. 
A depot is located at intersection node \#10 in the transportation network. Microgrid set comprises four microgrids located at bus \#14 in each distribution system and intersection nodes \#2, \#3, \#17, \#24 in the transportation network, respectively, as depicted in Fig. \ref{fig: integrated_test_system} and Fig. \ref{fig: 33_bus_system}. The parameters for microgrids and the MESS fleet are described in Table \ref{table: microgrids} and Table \ref{table: mess}, respectively. 

%

\begin{table}[!tbp]
	\centering
	\caption{Generation Resources and Local Loads for microgrids} \label{table: microgrids}
	\resizebox{0.9\columnwidth}{!}{
	\begin{tabular}{cccccc}
		\hline
		\multicolumn{2}{c}{Microgrid \#}             & 1 & 2 & 3 & 4 \\ \hline
		\multirow{4}{*}{Generation} & $\overline{P}_{\text{DG},m}$ (MW) & 1.80 & 1.60 & 1.80 & 1.60 \\ \cline{2-6}
		& $\overline{Q}_{\text{DG},m}$ (MVar) &  1.35 & 1.20 & 1.35 & 1.20 \\ \cline{2-6}
		& $\overline{E}_{\text{DG}, m}$ (MWh)     & 34.5 & 30.7 & 34.5 & 30.7  \\ \cline{2-6}
		& $\uline{E}_{\text{DG}, m}$ (MWh)  &   3.5 & 3.0 & 3.5 & 3.0   \\ \hline
		\multirow{3}{*}{Local load} & Peak load  (MW)     &   \multicolumn{4}{c}{0.5}   \\ \cline{2-6}
		& Power factor   &   \multicolumn{4}{c}{0.9}  \\ \cline{2-6}
		& Load type      & C & R  & I  & I \\  \hline
		
		\multicolumn{6}{l}{Note: C - commercial, R - residential, I - industrial}
	\end{tabular}
	}
\end{table}


\begin{table}[!tbp]
	\centering
	\vspace{-10pt}
	\caption{Parameters of the MESS fleet} \label{table: mess}
	\resizebox{1.03\columnwidth}{!}{
		\begin{tabular}{ccccccc}
			\hline
			\multicolumn{1}{c}{MESS \#}           & 
				\begin{tabular}[c]{@{}c@{}}$\overline{P}_\text{ch}^\omega$ 	\\	 /	$\overline{P}_\text{dch}^\omega$ \\ (MW) 	\end{tabular} & 
				\begin{tabular}[c]{@{}c@{}} $\overline{E}_\text{c}^\omega$   \\ (MWh) \end{tabular} &  
				\begin{tabular}[c]{@{}c@{}} Initial \\ SOC \\ (\%)  \end{tabular} &
				\begin{tabular}[c]{@{}c@{}} $\overline{SOC}^\omega$ \\  / $\uline{SOC}^\omega$ \\ (\%) \end{tabular}   &
				\begin{tabular}[c]{@{}c@{}} $\eta_\text{ch}^\omega / \eta_\text{dch}^\omega $ \\ (\%) \end{tabular}  &
				\begin{tabular}[c]{@{}c@{}} $V_\text{avg}$  \\ (km/h)   \end{tabular}     \\ \hline
		\multicolumn{1}{c}{1} & \multirow{3}{*}{0.5} & \multirow{3}{*}{2.0} & \multirow{3}{*}{50} & \multirow{3}{*}{90/10} 		
			&\multirow{3}{*}{95/95}
			& 20 \\ 
		\multicolumn{1}{c}{2} & 	&   &  &  &  & 30 \\
		\multicolumn{1}{c}{3} & 	&   &  &  &  & 40 \\ \hline

		\end{tabular}
	}
\end{table}

In the remaining section, three sub-cases are considered to show the effectiveness of MESS mobility for service restoration as follows.

Case \Romannum{1}-a) There are no MESSs;

Case \Romannum{1}-b) Allocation of MESSs;

Case \Romannum{1}-c) Dynamic scheduling of MESSs.

\vspace{6pt}
\noindent \textit{2) Simulation Results}
\vspace{6pt}

Table \ref{table: Comparison of 3 cases} compares the three cases in terms of objective value and load restoration.

\begin{table}[!tbp]
	\centering
	\caption{Comparison in Case \Romannum{1}} \label{table: Comparison of 3 cases}
	\resizebox{1.02\columnwidth}{!}{
	\begin{tabular}{ccccc}
		\hline
		\multicolumn{2}{c}{Results}    & Case \Romannum{1}-a) & Case \Romannum{1}-b) & Case \Romannum{1}-c) \\ \hline
		
		\multirow{5}{*}{\begin{tabular}[c]{@{}c@{}} Objective \\ values \\ (\$)\end{tabular}} 
					& Interruption cost   &     288535 &  272202 &  231662   \\ \cline{2-5} 
					& \begin{tabular}[c]{@{}l@{}}MG generation cost\end{tabular} & 58752  &  58752  &    58752     \\ \cline{2-5} 
					&  \begin{tabular}[c]{@{}l@{}}Transportation cost\end{tabular}  & 0  &  640  & 1920      \\ \cline{2-5}
					& \begin{tabular}[c]{@{}l@{}}Battery maintenance cost\end{tabular}  & 0  &  1022 &    2709         \\ \cline{2-5}
					& Total cost &  347287   &  332616 & 295043          \\ \hline

		\multirow{3}{*}{\begin{tabular}[c]{@{}c@{}}Load \\ restoration \\ (\%)\end{tabular}} & Critical     &   80.98     &   82.26   &   88.61    \\ \cline{2-5} 
		& Non-critical &  41.73    &  43.57     &   36.24    \\ \cline{2-5}
		& Total & 58.33 & 59.62 & 59.48 \\ \hline
	\end{tabular}
	}
\end{table}

\textit{Case \Romannum{1}-a) There are no MESSs}: The base case assumes that there are no MESSs, microgrids only use local generating resources for service restoration in distribution systems. The total cost is $\$347287$, restoration of critical, non-critical and total loads are 80.98\%, 41.73\%, 58.33\%, respectively.

\textit{Case \Romannum{1}-b) Allocation of MESSs}: In the second case,  resource allocation is introduced to dispatch MESSs from the depot to microgrids in the initial stage of restoration, then MESSs stay at the microgrid to coordinate with local stationary resources. The MESSs are dispatched to microgrids \#2, \#4 and \#2, respectively, to provide support for local resources. In comparison with Case \Romannum{1}-a), the total cost decreases by 4.22\% to \$332616 while incorporating MESSs only introduces about 1.4\% additional energy. The restoration of critical, non-critical and total loads are 82.26\%, 43.57\%, and 59.62\%, respectively. The results show the benefits of properly positioning MESSs to coordinate with microgrids for service restoration. However, after the allocation of MESSs from depots to microgrids, the MESSs is then served as stationary resources, the mobility and flexibility are underutilized.

\textit{Case \Romannum{1}-c) Dynamic scheduling of MESSs}: This case co-optimize the dynamic scheduling of MESSs, resource dispatching of microgrids and network reconfiguration. The total cost reduces by 11.29\% than Case \Romannum{1}-b) to \$295043.
All the microgrid generation costs in the three cases are the same, this is due to the fact that three cases have the same condition of energy capacity $\overline{E}_{\text{DG}, m}$ and minimum reserve $\uline{E}_{\text{DG}, m}$ in corresponding microgrids, and the results indicate that all the microgrids get to the minimum reserve and are fully utilized for service restoration.
The restoration of critical, non-critical and total loads are 88.61\%, 36.24\% and 59.48\%, respectively.  It is noted that the total load restoration is a little bit less than Case \Romannum{1}-c), this is because MESSs transport energy among microgrids and have more charging/discharging behaviors, thus having more charging/discharging losses.

\begin{figure}[!tbp]
	\centering
	\includegraphics[width=\columnwidth, clip]{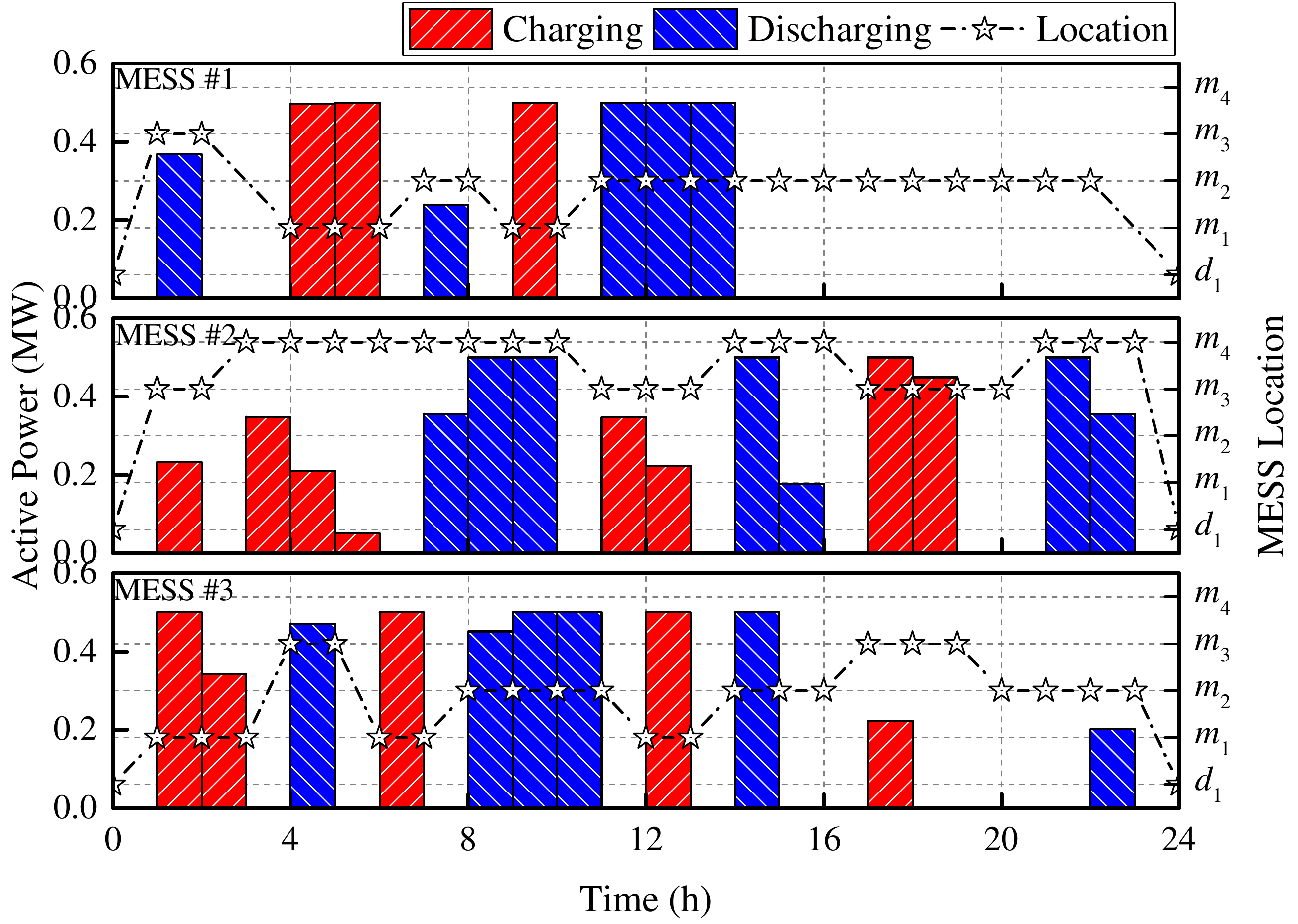}
	\caption{Scheduling results of the MESS fleet in Case \Romannum{1}-c).}
	\label{fig: scheduling_of_mess}
\end{figure}

Fig. \ref{fig: scheduling_of_mess} provides the charging/discharging schedule with respect to the position of MESS. The bar shows the charging/discharging active power while the dash lines with asterisks and right Y-axis indicates the MESSs' movements. The d1 represents the depot and m1-m4 represent microgrid \#1 - \#4, respectively. In contrast to the allocation of MESSs in Case \Romannum{1}-b), the dynamic scheduling of MESSs optimizes the sequence of movements and associated charging/discharging behaviors. Since the limited power and energy capacity, the MESSs need to move among microgrids back and forth considering the transportation cost.

For instance, MESS \#1 initially starts from depot \#1 to microgrid \#3  in (00:00-01:00) and discharges in (01:00-02:00). Then it takes three intervals to get to microgrid \#1 and stay there to charge in (04:00-06:00), followed by going to microgrid \#2 to discharge. Next, it moves back to microgrid \#2 to achieve energy transfer. Finally, it returns to depot \#1 in (22:00-24:00). It is also observed that the MESS \#1 is mainly dispatched back and forth between microgrid \#1 and \#2, MESS \#2 is mainly back and forth between microgrid \#3 and \#4, while MESS \#3 generally moves among more locations (microgrid \#1-\#3). This is due to the difference in the average speed of MESSs, as compared to MESS \#1 and \#2, MESS \#3 is faster and more effective to transfer energy among microgrids.

%
%
%

\begin{figure}[!tbp]
	\centering
	\includegraphics[width=\columnwidth, clip]{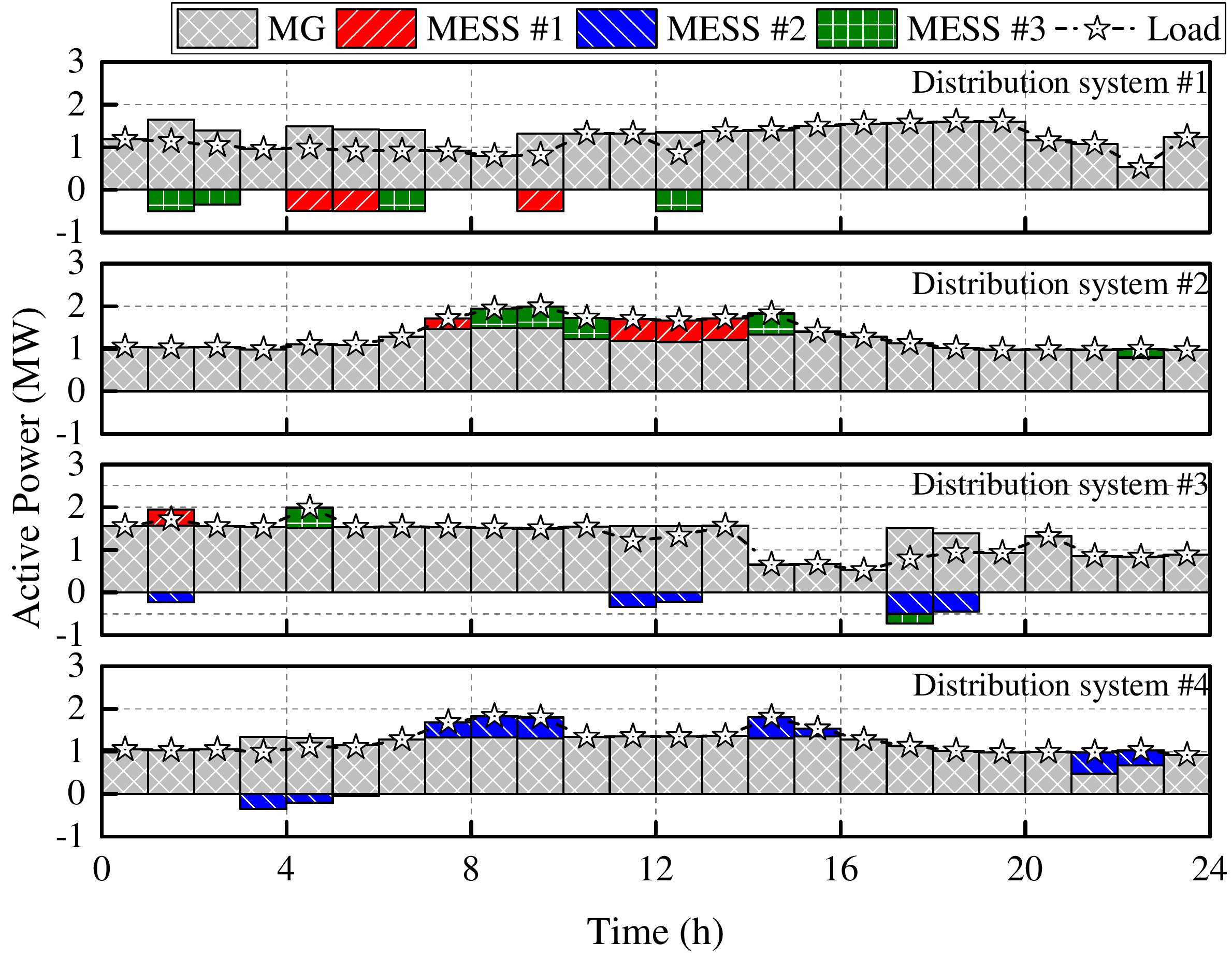}
	\caption{Coordinated generation dispatch and load restoration in Case \Romannum{1}-c).}
	\label{fig: generation_dispatch}
\end{figure}

\begin{figure}[!tbp]
	\centering
	\includegraphics[width=\columnwidth, clip]{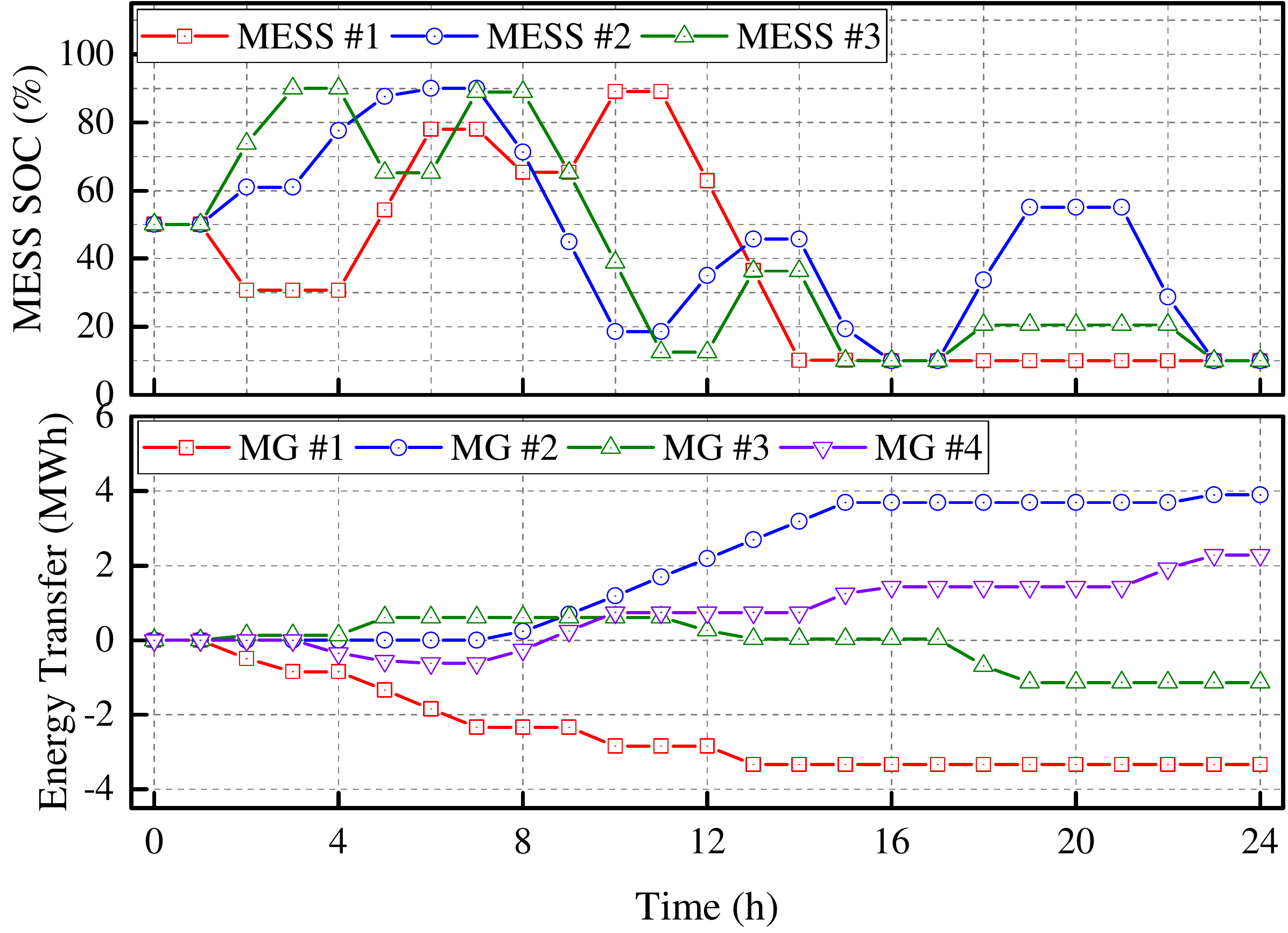}
	\caption{Energy transfer through MESSs in Case \Romannum{1}-c).}
	\label{fig: energy_transfer}
\end{figure}

Fig. \ref{fig: generation_dispatch} denotes the coordinated generation dispatch and load restoration in each distribution system. By coordinating multiple sources, the generation capacities are fully utilized to restore critical loads with higher customer interruption cost. Fig. \ref{fig: energy_transfer} describes the energy transfer among microgrids through charging or discharging of MESSs.  A microgrid with positive energy transfer means it receives energy from MESSs whereas negative one means outputting energy from this microgrid. It is observed that energy transfer is mainly from microgrids \#1, \#3 to microgrids \#2, \#4.

\vspace{4pt}
\noindent \textit{(a) Effect of Temporal-spatial Dynamics of MESSs}
\vspace{4pt}

The simulation demonstrates that MESSs supplement energy by charging from some microgrids with relatively surplus resources and transfer to other microgrids to minimize the total cost. During the restoration process, the energy imbalance between distribution systems is caused by topology and operation constraints, which impedes the effective utilization of local stationary resources. Therefore, the integration of MESSs and coordination with microgrids can make the most of MESSs mobility to deal with the energy imbalance. MESSs can transfer power and energy among microgrids to restore critical loads and reduces system total costs. 
It can be seen in Table \ref{table: Comparison of 3 cases} that total load restoration for Case \Romannum{1}-b) and Case \Romannum{1}-c) are almost the same whereas the load restoration for critical loads increases by 6.93\% and 5.37\% in Case \Romannum{1}-c) than Cases \Romannum{1}-a) and  \Romannum{1}-b) because MESSs mobility is better utilized to serve more critical loads. 
For instance, the MESS \#3 initially moves to microgrid \#1 from depot \#1 and charges in (01:00-03:00) to get to maximum SOC of 90\% (see Fig. \ref{fig: scheduling_of_mess} and Fig. \ref{fig: energy_transfer}). Next it moves to microgrid \#3 to discharge, with SOC reduced to 65.19\%. Then it moves back and forth between microgrid \#1 and \#2 in (06:00-16:00) to transfer energy from microgrid \#1 to \#2. The SOC ranges from 10\% to 88.94\% and MESS \#3 is fully discharged to the minimum SOC of 10\% at 15:00. Next it moves to microgrid \#3 to charge and back to microgrid \#2 to discharge. So the energy is mainly transferred from microgrids \#1 and \#3 to microgrid \#2. Similar temporal-spatial dynamics and associated charging/discharging behaviors can be observed for other MESSs.

%

Also, it is noticed that MESSs can perform load shifting within the same microgrid. For example, MESS \#2 charges at microgrid \#4 in (03:00-06:00) and discharges at the same location in (07:00-10:00). This is because the load is relatively low during some intervals, so MESS \#2 charges in energy sufficient hours and prepares for peak hours to obtain effective use of energy resources in microgrid \#4.

The comparison of three cases highlights the importance of effective utilization of MESSs mobility. The integrated restoration strategy reduces the total cost by coordinating the dynamic scheduling of MESSs, resource dispatching of microgrids and distribution network reconfiguration. In comparison with the allocation of MESSs only in the very initial stage of the restoration process, the MESSs mobility is fully utilized by dynamic scheduling to deal with energy imbalance in distribution systems posed by topology and operation constraints.


\vspace{4pt}
\noindent \textit{(b) Impact of Damage and Repair to Branches}
\vspace{4pt}

\begin{figure}[!tbp]
	\centering
	\includegraphics[width=\columnwidth, clip]{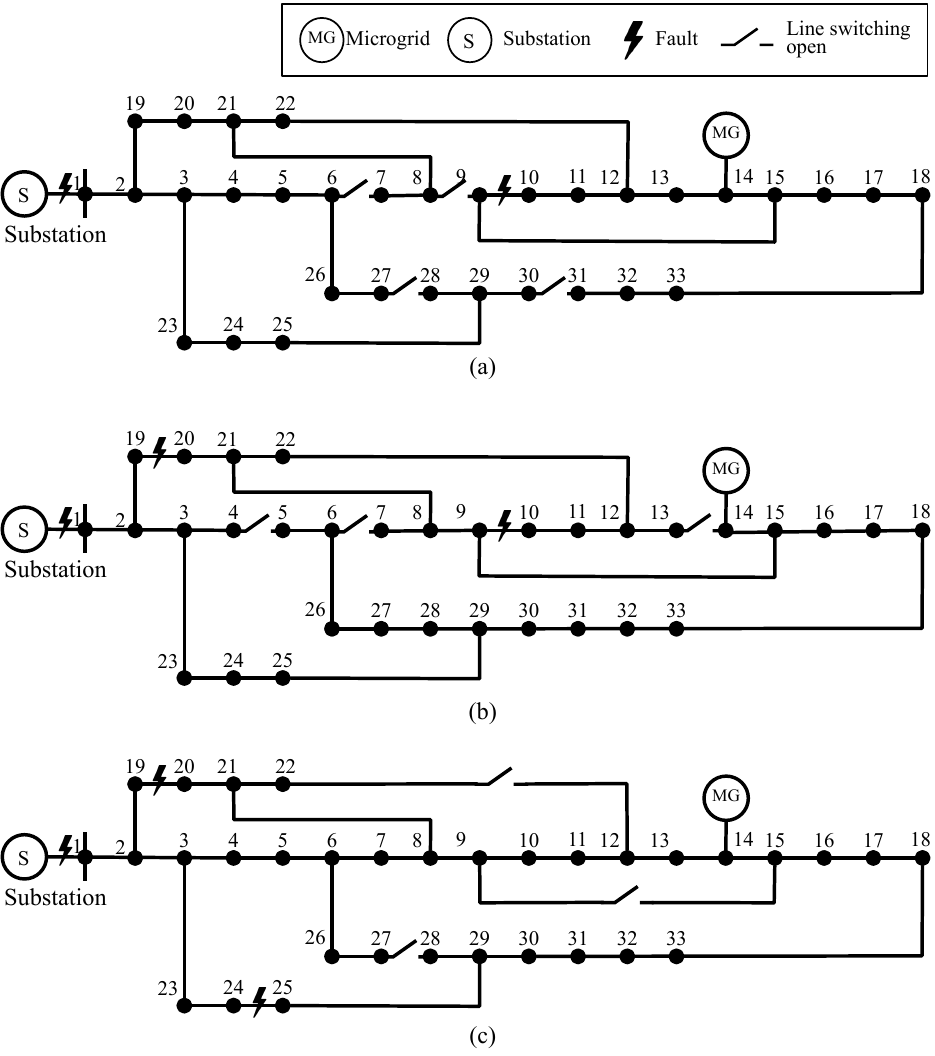}
	\caption{Network reconfiguration for distribution system \#3: (a) t=5, faults: substation, branch (9, 10); (b) t=6, faults: substation, branches (9, 10), (19, 20); (c) t=10, faults: substation, branches (19, 20), (24, 25), repair: branch (9, 10).}
	\label{fig: network_reconfiguration}
\end{figure}

To consider the subsequent damage and repair to branches in distribution systems, the system damage status is updated at each interval and the network topology is reconfigured, as shown in Fig. \ref{fig: network_reconfiguration}. 
For example, there are two faults in distribution system \#3 at $t=5$: substation fault and damage to branch (9, 10). The opening line switches are (6, 7), (8, 9), (27, 28), (30, 31), as shown in Fig. \ref{fig: network_reconfiguration}(a). At $t=6$, new damage to branch (19, 20) occurs, and the distribution system \#3 is reconfigured to new topology with opening switches (4, 5), (6, 7), (13, 14), as shown in Fig. \ref{fig: network_reconfiguration}(b). Similarly, subsequent damage to branches (24, 25), (18, 33) change the network topology and reconfiguration is needed. 

In addition, the repair to branches is considered. Once the repair is finished, the branch can be controlled again. For example, as shown in Fig. \ref{fig: network_reconfiguration}(c), branch (9, 10) is repaired at $t=10$, and the topology is reconfigured by opening switches (9, 15), (12, 22), (27, 28). At $t=20$, branch (24, 25) is repaired, the distribution network takes the new topology and switches (6, 26), (8, 21), (9, 10), (12, 22) are open.


\vspace{4pt}
\noindent \textit{(c) Impact of Damage and Repair to Roads}
\vspace{4pt}

\begin{figure}[!tbp]
	\centering
	\includegraphics[width=\columnwidth, clip]{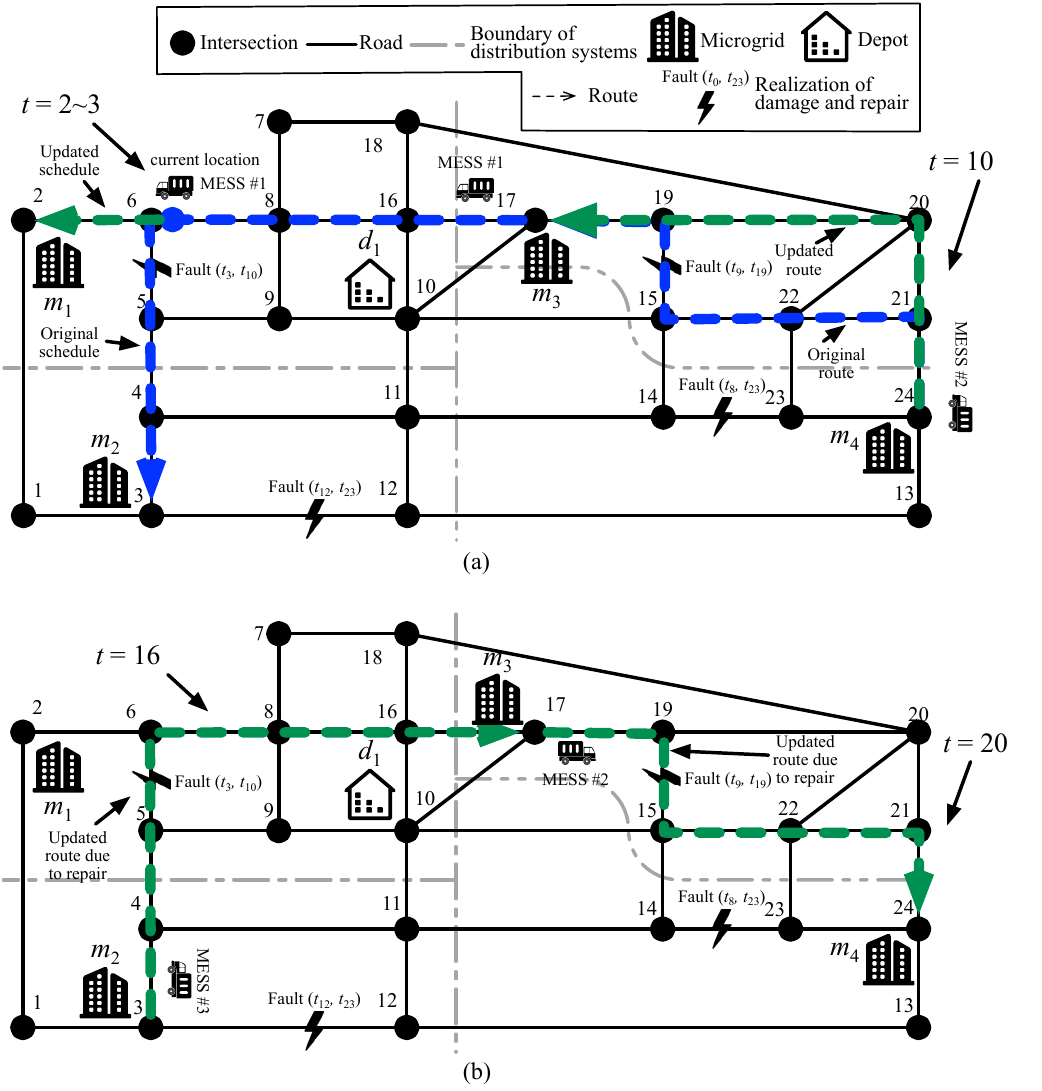}
	\caption{Impact analysis of damage and repair to roads: (a) impacts of damage to roads on vehicle scheduling and routing, (b) impacts of repair to roads on vehicle scheduling and routing.}
	\label{fig: damage_and_repair_roads}
\end{figure}

Extreme events also cause damage to transportation networks, which in turn will impact the scheduling of the MESS fleet. 
In order to leverage the updated information on subsequent damage and repair to roads,  the scheduling of the MESS fleet is optimized at each interval over the prediction horizon based on the updated information and current locations. The scheduling results are shown in Fig. \ref{fig: damage_and_repair_roads}(a). In general, MESSs have two choices when damage occurs, one is to discard the original scheduling and be rescheduled to another microgrid, the other option is to be rerouted to the original destination via a different route. 
For example, the MESS \#1 is dispatched from microgrid \#3 to microgrid \#2 at $t=2$, as depicted by the blue dot line. At $t=3$, the current location of MESS \#1 is obtained, shown as the blue node in Fig. \ref{fig: damage_and_repair_roads}(a). The subsequent damage causes road 5-6 to fail, thus the originally designated route is no longer available. The optimization takes into account the updated information and the current location, consequently MESS \#1 is rescheduled to microgrid \#1 via a new route, as depicted by the green dot line. 
Meanwhile, MESSs can be rerouted to the same destination when damage occurs. The originally designated route for MESS \#2 from microgrid \#4 to microgrid \#3 is shown by the blue dash line. The road (15, 19) is damaged at $t=10$, so the MESS \#2 takes another route to microgrid \#3, which is depicted by the green dash line in Fig. \ref{fig: damage_and_repair_roads}(a). 

In addition, Fig. \ref{fig: damage_and_repair_roads}(b) illustrates the impact of repair to roads. Once a road is repaired, scheduling and routing of the MESS fleet can be updated. At $t=16$, the road 5-6 is repaired and available again, thus the MESS \#3 is dispatched from microgrid \#2 to \#3 via the designated route, which is represented by the green dash line. Similarly, after the repair of road 15-19, MESS \#2 is dispatched from microgrid \#3 to microgrid \#4 via the route with the shortest travel time.

The results reveal the flexibility of the proposed model that MESS fleets can be rescheduled or rerouted  at each interval considering the damage and repair to roads.

\subsection{Case \Romannum{2}: Singapore Transportation Network with Six 33-bus Distribution Systems}

\vspace{6pt}
\noindent \textit{1) Test Systems}
\vspace{6pt}

\begin{figure}[!tbp]
	\centering
	\includegraphics[width=\columnwidth, clip]{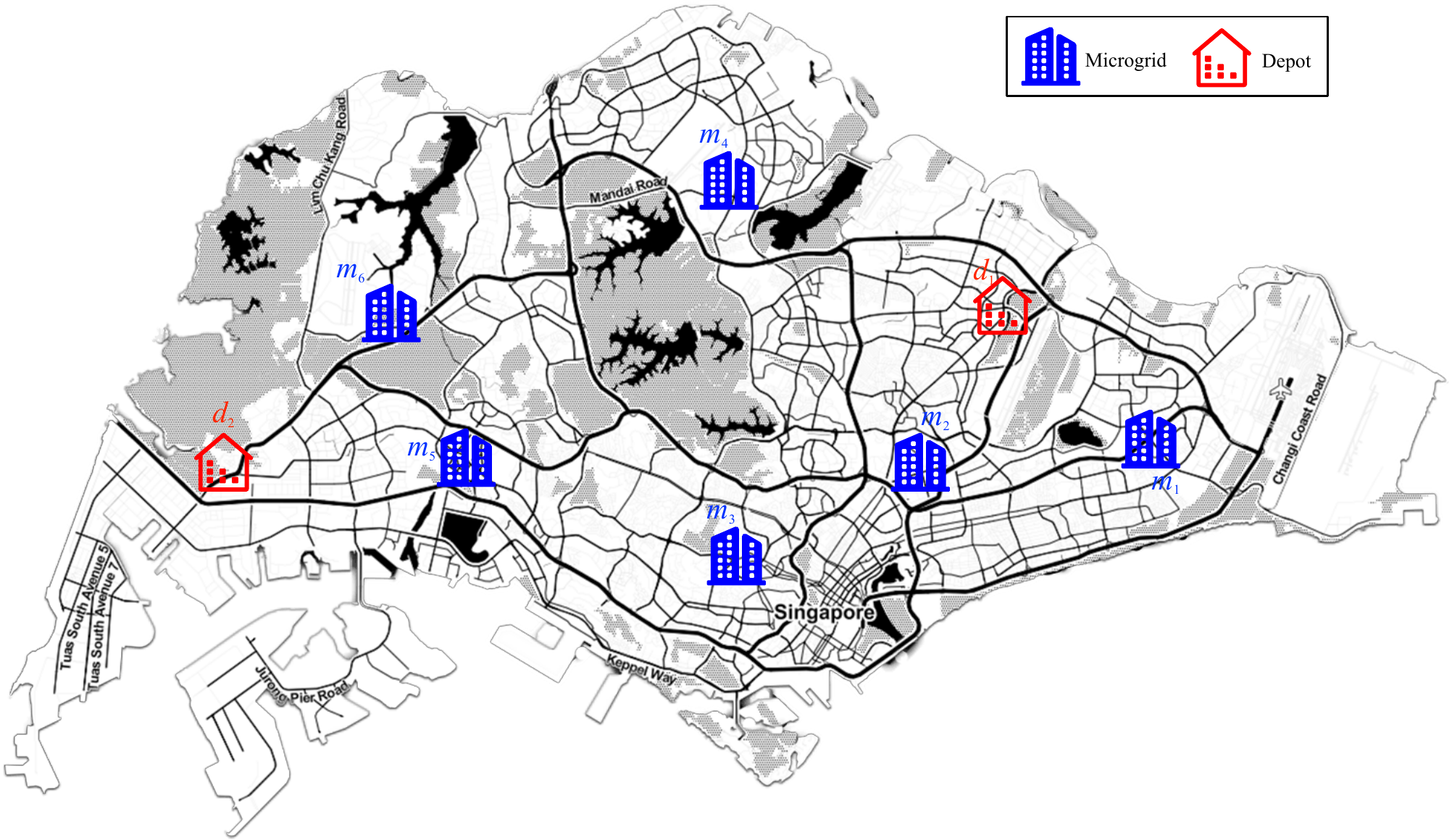}
	\caption{Singapore transportation network with microgrids and depots.}
	\label{fig: singapore_network}
\end{figure}

To verify the scalability of the proposed integrated restoration strategy, the case study on an integrated test system with Singapore transportation network connecting six 33-bus distribution systems is carried out for illustration. Each distribution system has one microgrid located at bus 14. Fig. \ref{fig: singapore_network} shows the transportation network with microgrids and depots' locations.
The Python client for Google Maps API \cite{Google} is used to retrieve geospatial data for Singapore transportation network. Six microgrids and two depots are located across the map. Microgrids \#5 and \#6 have the same properties as Microgrids \#1 and \#3, respectively. A fleet of five MESSs is considered, with three MESSs initially located at depot \#1, and the other two at depot \#2. Other parameter settings follow the Case \Romannum{1}-c). Similarly, three cases are implemented as follows.

Case \Romannum{2}-a) There are no MESSs;

Case \Romannum{2}-b) Allocation of MESSs;

Case \Romannum{2}-c) Dynamic scheduling of MESSs.

\vspace{6pt}
\noindent \textit{2) Simulation Results}
\vspace{6pt}

\begin{table}[!tbp]
	\centering
	\caption{Comparison in Case \Romannum{2}} \label{table: Comparison of 3 cases_large}
	\resizebox{1.02\columnwidth}{!}{
		\begin{tabular}{ccccc}
			\hline
			\multicolumn{2}{c}{Results}    & Case \Romannum{2}-a) & Case \Romannum{2}-b) & Case \Romannum{2}-c) \\ \hline
			
			\multirow{5}{*}{\begin{tabular}[c]{@{}c@{}} Objective \\ values \\ (\$)\end{tabular}} 
			& Interruption cost   &     362892 &  344747 & 280455      \\ \cline{2-5} 
			& \begin{tabular}[c]{@{}l@{}}MG generation cost\end{tabular} &   87990  & 87990 &    87990      \\ \cline{2-5} 
			&  \begin{tabular}[c]{@{}l@{}}Transportation cost\end{tabular}  &  0   &  800 & 2320       \\ \cline{2-5}
			& \begin{tabular}[c]{@{}l@{}}Battery maintenance cost\end{tabular}  &   0     &  1694    &     3947   \\ \cline{2-5}
			& Total cost &   450882   &  435231    &    374712      \\ \hline
			
			\multirow{3}{*}{\begin{tabular}[c]{@{}c@{}}Load \\ restoration \\ (\%)\end{tabular}} & Critical     &   81.12    &    	85.62    &  91.64     \\ \cline{2-5} 
			& Non-critical &   51.03   &    56.48   &   48.92    \\ \cline{2-5}
			& Total & 60.01 &  61.85 &  61.66 \\ \hline
		\end{tabular}
	}
\end{table}

\begin{figure}[!tbp]
	\centering
	\includegraphics[width=\columnwidth, clip]{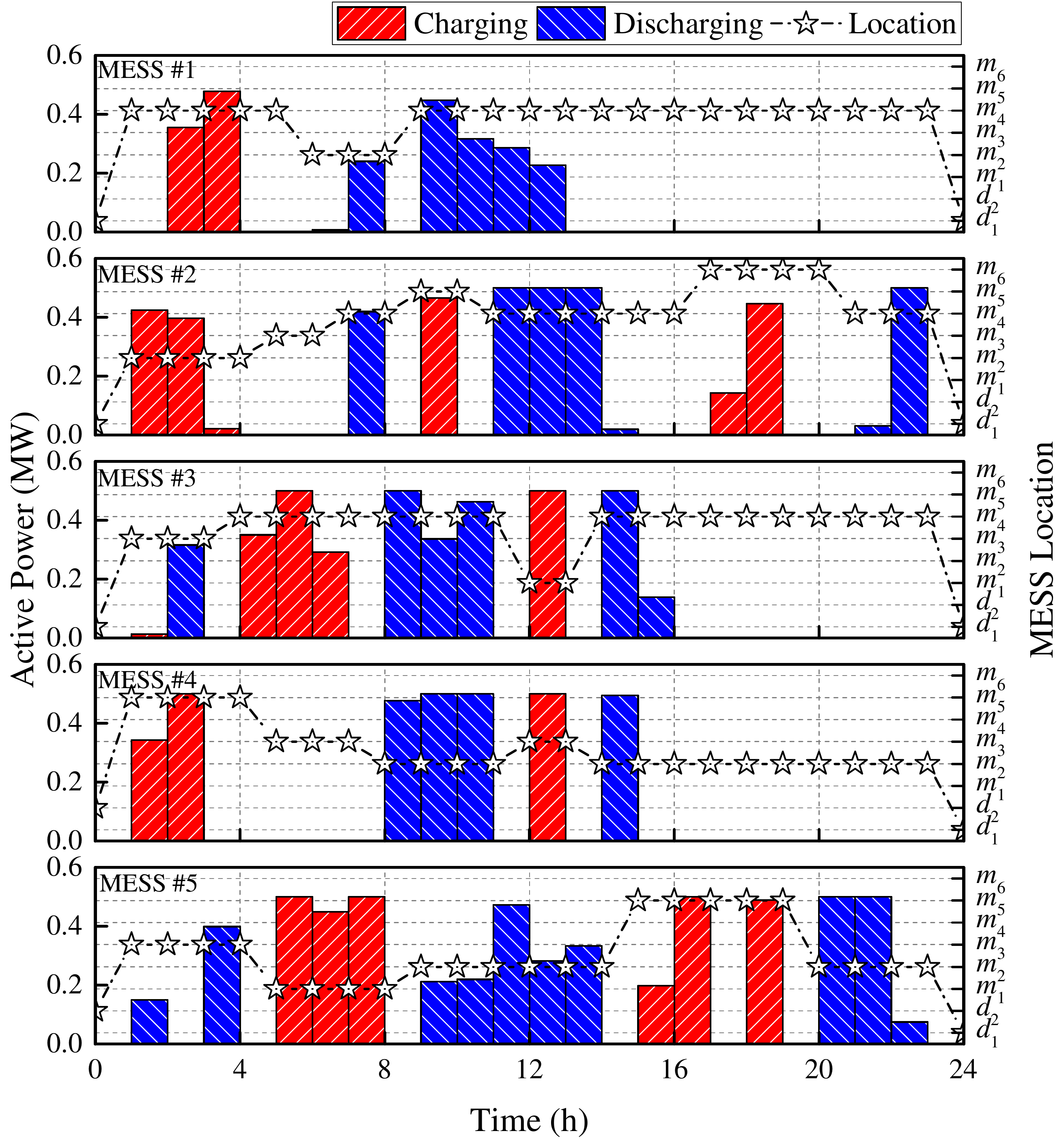}
	\caption{Scheduling results of the MESS fleet in Case \Romannum{2}-c).}
	\label{fig: scheduling_of_mess_large}
\end{figure}

\begin{figure}[!tbp]
	\centering
	\includegraphics[width=\columnwidth, clip]{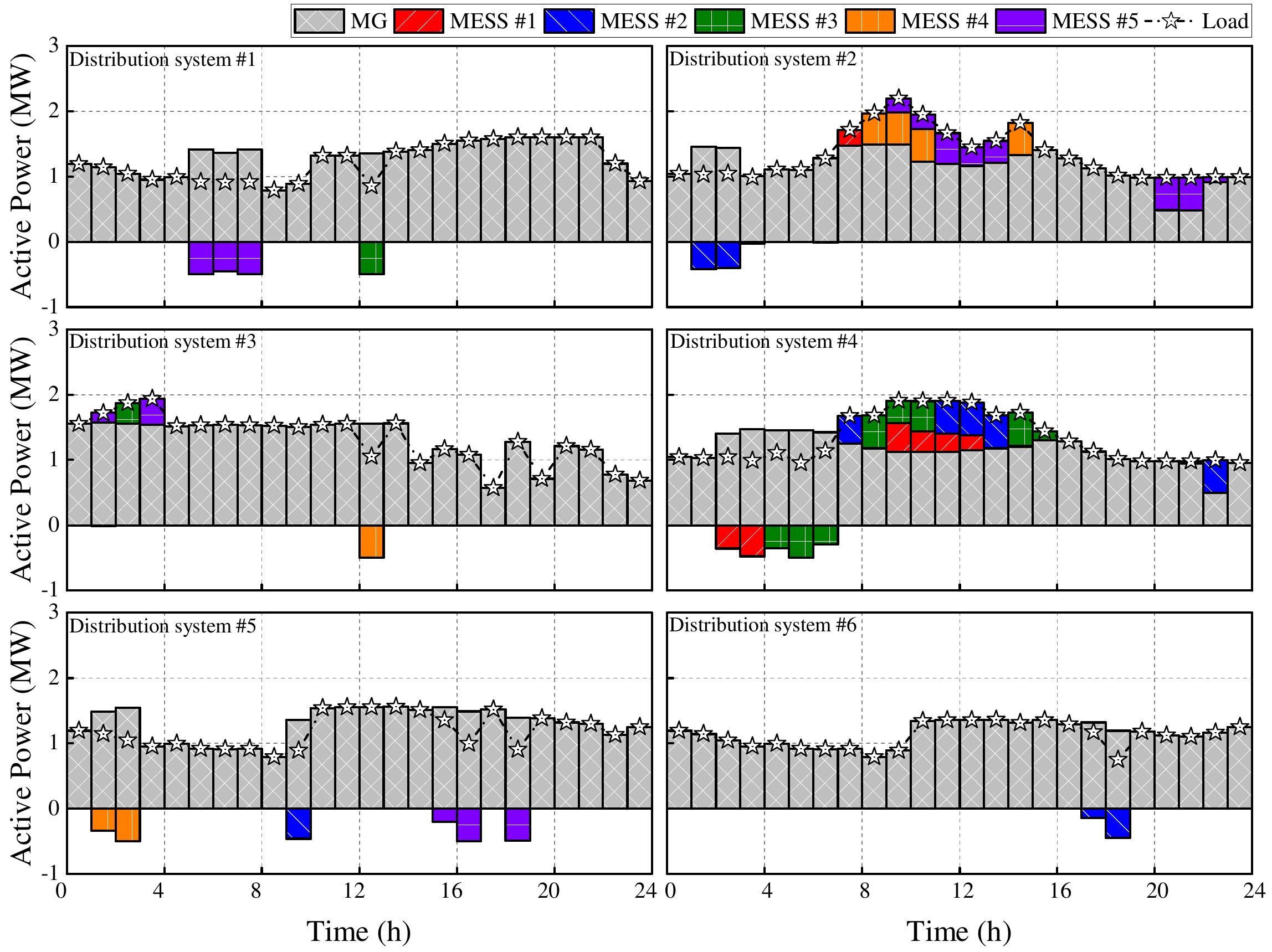}
	\caption{Coordinated generation dispatch and load restoration in Case \Romannum{2}-c).}
	\label{fig: generation_dispatch_large}
\end{figure}

Table \ref{table: Comparison of 3 cases_large} depicts the objective value and load restoration in Case \Romannum{2}. 
The total costs in Case \Romannum{2}-c) \ reduce by 16.89\% than Case \Romannum{2}-a) and 13.91\% than Case \Romannum{2}-b). In Case \Romannum{2}-c), the restoration of critical, non-critical and total loads are 91.64\%, 48.92\% and 61.66\%, respectively. More critical loads with higher importance are restored.
Fig. \ref{fig: scheduling_of_mess_large} depicts the scheduling results for the MESS fleet. The dynamic scheduling optimizes the sequence of movements and charging/discharging behaviors. MESSs get charged from some microgrids with surplus resources and transport energy to other microgrids to minimize the total system cost. 
Fig. \ref{fig: generation_dispatch_large} shows the coordinated generation dispatch and load restoration. By coordinating stationary and mobile resources, the generation capacities of microgrids and MESSs are better utilized for service restoration.


The simulation results demonstrate the potential applications of the proposed integrated restoration strategy to enhance distribution system resilience.

\section{Conclusions}
\label{sec: Conclusions}

This paper proposes a rolling integrated service restoration strategy to minimize the total system cost by coordinating MESS fleets, microgrids and distribution systems. The proposed service restoration strategy takes into account damage and repair to both the roads in transportation systems and the branches in distribution systems. 
The uncertainties in load consumption and the status of roads and branches are considered to generate scenarios using Monte Carlo simulation method. 
A rolling optimization framework is adopted to consider subsequent damage and repair during the restoration process. The operation of MESS fleets is modeled by a stochastic multi-layer time-space network technique and MESS fleets can be rescheduled and rerouted at each interval.
The coordinated scheduling at each interval over the prediction horizon is formulated as a two-stage stochastice MILP. 
The simulation results demonstrate the effectiveness of MESSs mobility that transfers energy across multiple distribution systems and coordinates with microgrids during the disasters, and highlight that the mobile and stationary resources can be well coordinated to enhance distribution system resilience. 

%
%
%
%
%



\ifCLASSOPTIONcaptionsoff
  \newpage
\fi



%

\normalem

\bibliographystyle{IEEEtran}
\bibliography{rotess_ref.bib}

%

\begin{IEEEbiography}[{\includegraphics[width=1in,height=1.25in,clip,keepaspectratio]{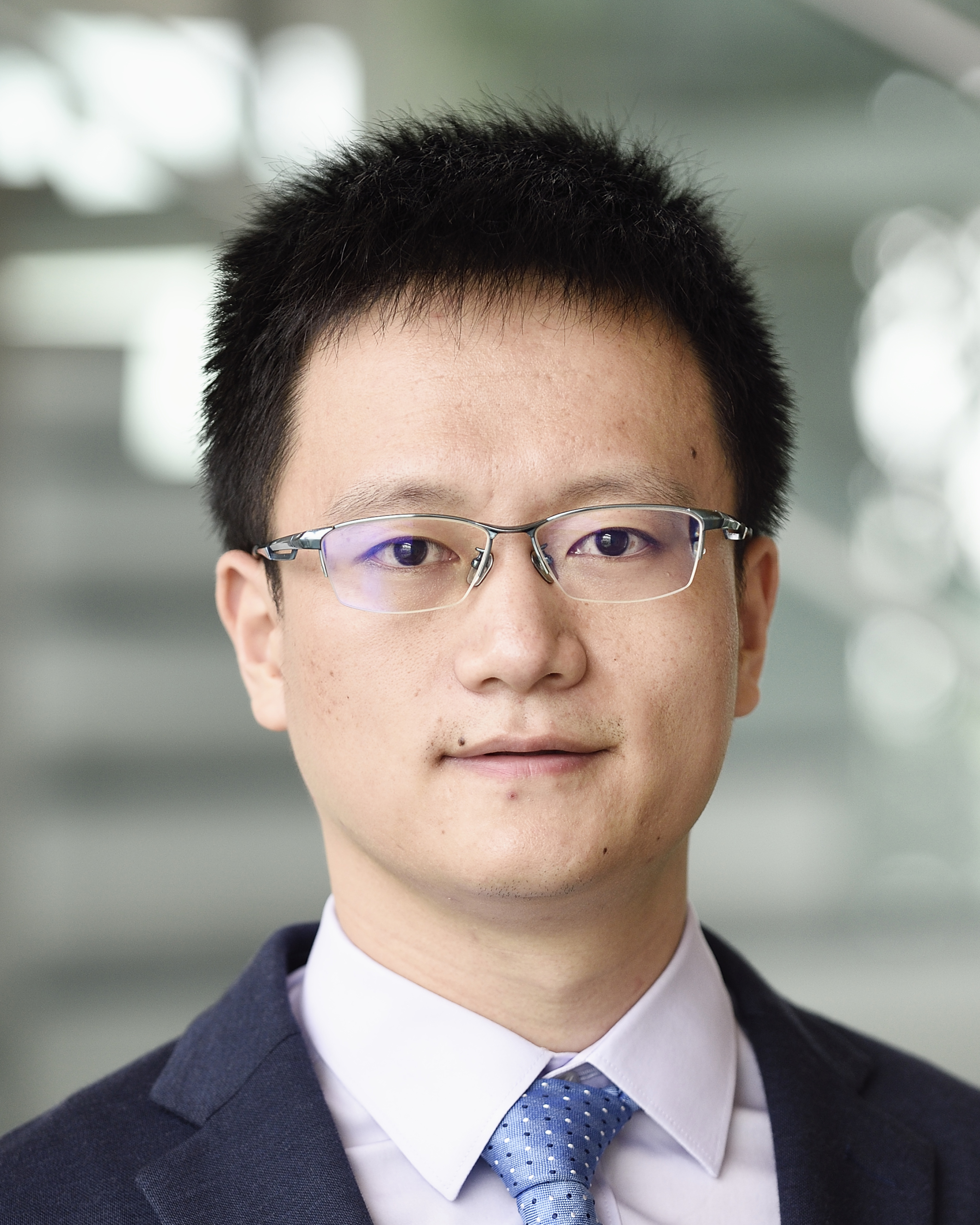}}]
{Shuhan Yao} (S'16) received the B.Eng. degree in electrical engineering from Chongqing University, Chongqing, China, in 2013. Currently, he is pursuing the Ph.D. degree in the Interdisciplinary Graduate School, Nanyang Technological University, Singapore. His research interests include smart grid resilience, optimization and reinforcement learning.
\end{IEEEbiography}

\vfill

\begin{IEEEbiography}[{\includegraphics[width=1in,height=1.25in,clip,keepaspectratio]{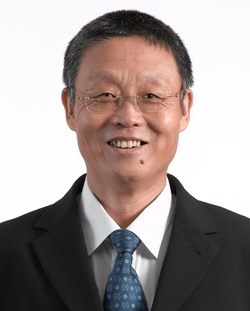}}]
{Peng Wang} (F’18) received the B.Sc. degree in electronic engineering from Xian Jiaotong University, Xian, China, in 1978, the M.Sc. degree from Taiyuan University of Technology, Taiyuan, China, in 1987, and the M.Sc. and Ph.D. degrees in electrical engineering from the University of Saskatchewan, Saskatoon, SK, Canada, in 1995 and 1998, respectively. Currently, he is a Professor with the School of Electrical and Electronic Engineering at Nanyang Technological University, Singapore.
\end{IEEEbiography}

\begin{IEEEbiography}[{\includegraphics[width=1in,height=1.25in,clip,keepaspectratio]{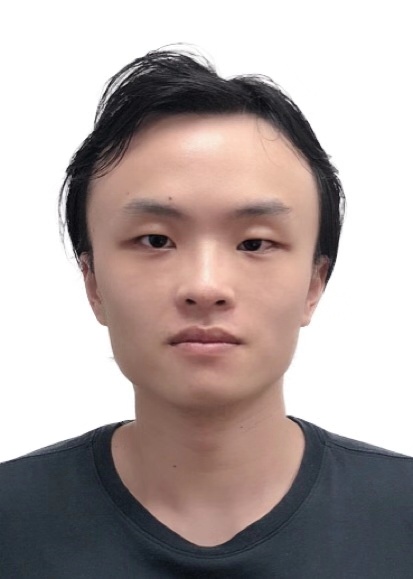}}]
{Xiaochuan Liu} (S’17) received the B.Eng. degree in electrical engineering from Wuhan University, Wuhan, China, in 2015, and the M.Sc. degree from Nanyang Technological University, Singapore, in 2016. Currently, he is pursuing his Ph.D. degree in the School of Electrical and Electronic Engineering, Nanyang Technological University, Singapore. His research interests include applied optimization and power system operation.
\end{IEEEbiography}

\begin{IEEEbiography}[{\includegraphics[width=1in,height=1.25in,clip,keepaspectratio]{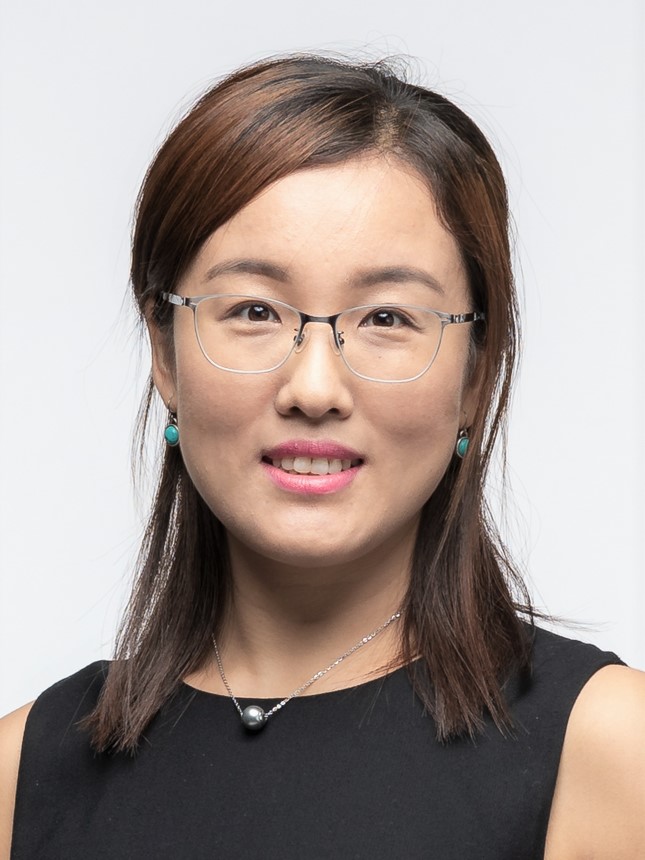}}]
{Huajun Zhang} (S’18) received the B.Eng., M.Eng. degree in electrical engineering from Shandong University, Jinan, China, in 2008 and 2011, respectively. Currently, she is pursuing the Ph.D. degree in the Interdisciplinary Graduate School, Nanyang Technological University, Singapore. Her research interest includes reliability and resilience evaluation of multi-energy systems.
\end{IEEEbiography}

\begin{IEEEbiography}[{\includegraphics[width=1in,height=1.25in,clip,keepaspectratio]{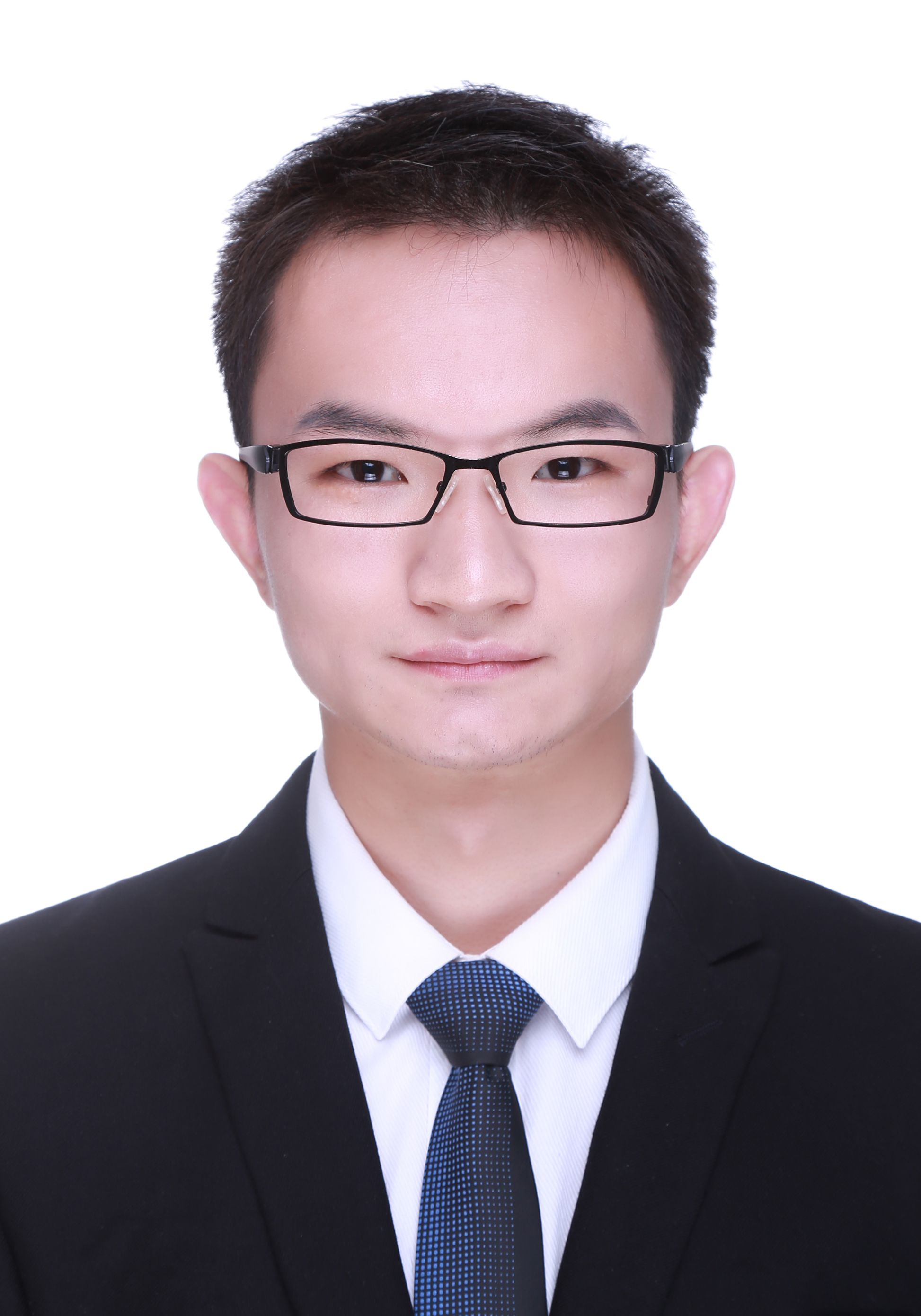}}]
{Tianyang Zhao} (M’18) received the B.Sc., M.Sc. and Ph.D. degree in automation of electric power systems from North China Electric Power University, Beijing China in 2011, 2013 and 2017, respectively. Currently, he is post-doc research fellow at Energy Research Institute @ Nanyang Technological University (ERI@N), Singapore. His research interest includes power system operation optimization and game theory.
\end{IEEEbiography}

\vfill


\end{document}